\documentclass[12pt,reqno]{amsproc}
\usepackage[centertags]{amsmath}
\usepackage{amssymb,amsfonts,amsthm}
\newcommand{\TryPackage}[3]{\IfFileExists{#1.sty}
{\usepackage{#1}#2}{#3}}\TryPackage{mathrsfs}
{\renewcommand{\mathcal}{\mathscr}}
{
\TryPackage{eucal}{}{}}
\usepackage{a4mlnewusa,math40,local}

\begin{document}

\title[On the number of square integrable solutions ...]{On the number of square integrable solutions
and self--adjointness of symmetric first order
systems of differential equations}

\author{Matthias Lesch}
\address{The University of Arizona\\
Department of Mathematics\\ 
617 N. Santa Rita\\
Tucson, AZ, 85721--0089\\
USA}
\email{lesch@math.arizona.edu}
\thanks{The first named author was supported by a Heisenberg
fellowship of Deutsche Forschungsgemeinschaft and by the
National Science Foundation under Grant No.
DMS--0072551}

\author{Mark Malamud}
\address{Department of Mathematics\\ University of 
Donetsk\\ Donetsk\\ 
Ukraine}
\email{mmm@univ.donetsk.ua}

\begin{abstract}
The main purpose of this paper is to investigate the formal 
deficiency indices $\cn_{\pm}(I)$ of a symmetric first order system
$$
Jf'+Bf=\gl \ch f
$$
on an interval $I$, where $I=\R$ or $I=\R_\pm.$  Here $J,B,\ch$ 
are $n\times n$ matrix valued functions and the Hamiltonian $\ch\ge 0$ may be
singular even everywhere. We obtain two results for such a system to have 
minimal numbers $\cn_\pm(\R)=0$ (resp. $\cn_\pm(\R_\pm)=n$)
and a criterion for their maximality 
$\cn_{\pm}(\R_+)=2n.$ Some conditions for a canonical system to have 
intermediate numbers $\cn_\pm(\R_+)$ are presented, too. We also obtain 
a generalization of the well--known Titchmarsh--Sears 
theorem for second order
Sturm--Liouville type equations. This contains results due to Lidskii and Krein
as special cases.

We present two approaches to the above problems: one dealing with
formal deficiency indices and one dealing with (ordinary) deficiency
indices. Our main (non--formal) approach is based on the investigation of a 
symmetric linear 
relation $\Smin$ which is naturally associated to a first order system. 
This approach works in the framework of extension theory and therefore
we investigate in detail
the domain $\cd(\Smin^*)$ of $\Smin^*$. In particular, we prove the so called 
regularity theorem for $\cd(\Smin^*)$.

The regularity result allows us to construct a bridge between the
"formal" and "non--formal" approaches by establishing a connection between the
formal deficiency indices $\cn_{\pm}$ and the usual deficiency
indices $N_{\pm}(\Smin)$. In particular we have 
$\cn_{\pm}=N_{\pm}$ for definite systems.


As a byproduct of the the regularity result we obtain
very short proofs of (generalizations of) the main 
results of the paper by Kogan and Rofe--Beketov \cite{KogRof:SIS} as well as a 
criterion for the quasi--regularity of canonical systems. This covers the 
Kac--Krein theorem and some results from \cite{KogRof:SIS}.
\end{abstract}

\maketitle
\tableofcontents

\section{Introduction}\label{sec1}
Let $I\subset\R$ be an interval and consider the first order system
      \begin{equation}\label{intro-G1}
J(x)f'(x)+B(x)f(x)=\ch(x)g(x),
     \end{equation}
where $J,B,\ch:I\rightarrow \Mat(n,\C)$ are locally integrable matrix--valued
functions with $\ch\ge 0$ and $J(x)$ invertible
(cf. \eqref{G1.0} below for the precise assumptions on $J,B,\ch$).

We would like to consider $g$ in \eqref{intro-G1} as the result obtained
by applying an operator to $f$. However, certain difficulties arise
if $\ch\ge 0$ is singular. It turns out that the appropriate framework 
to study \eqref{intro-G1} is the framework of symmetric linear relations 
in Hilbert space (Def. \plref{SDM-S1.1}). To outline this let $\cLH$ be 
the space of $\C^n$--valued measurable
functions $f$ with $\int_I f^*\ch f<\infty$ and denote by $\LH$
the corresponding Hilbert space (equivalence classes!). Then \eqref{intro-G1}
induces symmetric linear relations, $\cs, S$, in the spaces
$\cLH,\LH$ in a natural way. The first major problem which arises is the 
regularity problem.
Suppose that one has classes $\tilde f,\tilde g\in \LH$ such that 
$\{\tilde f,\tilde g\}\in S$. Are there representatives $f,g\in\cLH$ of
$\tilde f,\tilde g$ such that \eqref{intro-G1} holds? In this case $f$ would
be automatically absolute continuous, because $J(x)$ is
invertible. Therefore, it is appropriate to address this problem as regularity
problem. 

We answer this problem affirmatively (Theorem \plref{regthm}), generalizing 
work of Orcutt \cite[Thm. II.2.6 and 
Thm. IV.2.5]{Orc:CDE} and I. S. Kac \cite{Kac:LRGC},  \cite{Kac:LRG}.

The other major purpose of this paper is to generalize several criteria for
essential self--adjointness of first and second order differential operators 
to the present setting. 
We present two approaches to the above problems: one dealing with
formal deficiency indices and one dealing with (ordinary) deficiency
indices. Our main (non--formal) approach is based on the investigation 
of a symmetric linear 
relation $\Smin$ which is naturally associated to a first order system. 
This approach works in the framework of extension theory and therefore
we investigate in detail
the domain $\cd(\Smin^*)$ of $\Smin^*$. In particular, we prove the so called 
regularity theorem for $\cd(\Smin^*)$.

More precisely, the paper is organized as follows:

In Section \ref{sec2} we give a brief overview of the theory of symmetric
first order systems and introduce symmetric linear relations associated with
such a system. We present examples which show that on the one hand such s.l.r. may
have a very exotic behavior (Example \plref{ML-S2.2}) and on the other hand
that they occur quite naturally (Example \plref{SDM-S1.2}). Moreover,
we state the regularity Theorem (Theorem \plref{regthm}) and discuss various
normal forms of symmetric first order systems using gauge transformations. 
For the latter we follow Kogan and Rofe--Beketov \cite{KogRof:SIS}.

In Subsection \ref{sec21} we investigate the properties of $\cs, S$ on finite
intervals. The results, in particular regularity, are summarized in 
Proposition \plref{S1.6}.
The case of an infinite interval is presented in Subsection \ref{sec22}
(Proposition \plref{S1.8}).

So called \emph{definite} systems have more pleasant properties than general
systems. In Subsection \plref{sec23} we briefly discuss such systems and
present a criterion for definiteness.

Subsections \plref{sec24} and \plref{sec25} are devoted to defect spaces and
deficiency indices. Analogously as for a symmetric operator in a Hilbert space
the deficiency indices of a symmetric linear relation determine whether
it is essentially self--adjoint resp. whether there exist self--adjoint
extensions. In the case of the relations $\cs$ and $S$ associated to a first
order system one has to distinguish between the deficiency indices 
$N_{\pm}(S)$ of the s.l.r. $S$ in the Hilbert space 
$\LH$ and the \emph{formal deficiency indices}
$\cn_{\pm}(\cs)$ of the relation $\cs$ in the linear space $\cLH$. 
The latter is the dimension of \emph{formal defect subspace}
$\ce_\gl:=\setdef{f\in\cLH}{Jf'+Bf=\gl\ch f}, \gl=\pm i.$ 

For arbitrary systems we establish (Proposition \plref{ML-S2.19}) 
the equalities  $\dim\ce_{\pm\gl}(S)=N_{\pm}(S) + n -\rank S$ which turn 
into the equalities $N_{\pm}=\dim\ce_{\pm\gl},\gl \in \Bbb C_{\pm}$, 
for definite systems.     

This yields in particular that $\dim\ce_{\gl}$ is locally
constant in $\C\setminus\R$ for an arbitrary (not necessarily definite) 
system on an arbitrary interval (Proposition \plref{S1.11}). 
For $I=\R_+$ and, under more restrictive
assumptions for $I=\R$, this fact is due to Kogan
and Rofe--Beketov \cite[Theorem 2.1, Theorem 2.3]{KogRof:SIS}. 
On the one hand Proposition \plref{S1.11} improves    
\cite[Theorem 2.3]{KogRof:SIS}    
and on the other hand it gives a new proof of 
\cite[Theorem 2.1]{KogRof:SIS}                   
which is considerably simpler than the original proof. 
Our proof depends, however, on the regularity Theorem \plref{regthm}.

In Section \plref{sec3} we discuss essential self--adjointness of
the s.l.r. $S$ on the line. The essential self--adjointness criterion
Theorem \plref{S2.2} requires that $\ch$ is positive definite on
a sufficiently large set. 
In Subsection \plref{sec32} we deal with
the case in which \eqref{intro-G1} defines a symmetric operator.

The supplementary Section \plref{sec4} is included for completeness. 
We present an alternative proof of Theorem \plref{S2.2} using the 
well--known hyperbolic equation method.

Finally, Section \plref{sec5} discusses in more detail the deficiency 
indices of the system $S$ on the half--line. 
Here using simple arguments based on J. von Neumann formula 
we establish a connection between deficiency indices of the system 
$S$ considered on the half--lines $\Bbb R_{\pm}$ and on the 
line respectively. 

Combining this formula with the regularity results from Section 2
one immediately obtains the corresponding formula for the formal
deficiency indices from \cite[Sec. 2.3]{KogRof:SIS}.
Moreover, we generalize \cite[Sec. 2.3]{KogRof:SIS}
since our formula holds for arbitrary (not necessarily definite)
systems. This formula allows to translate results on
the half--line (about (formal) deficiency indices) into corresponding
results for the line and vice versa. In particular Theorem \plref{th3.1}
corresponds to Theorem \plref{S2.2}. However, in Subsection \plref{sec5.1}
we present a proof independent of Theorem \plref{S2.2}.


In Subsection \plref{sec5.2} we present a criterion for essential 
self--adjointness in a case where the Hamiltonian $\ch$ is singular 
(Theorem \plref{ML-S3.2.2}). This applies in particular
to second order Sturm--Liouville type equations. Our criterion generalizes
result's due to Lidskii \cite{Lid:NSI} and Krein \cite{Kre:TFO}
and it is in the spirit of the well--known
Titchmarsh--Sears theorem \cite{BerShu:SE}.


Furthermore, in Subsections \plref{sec3.3} and \plref{sec5.4} 
we present several other criteria which allow to determine the deficiency 
indices on the half line in several cases. In particular, 
Theorem \plref{th3.2} and Corollary \plref{cor5.18A}    
state a necessary and sufficient condition for
a first order system to have maximal deficiency indices as well as to be
quasiregular. These criteria have been inspired by the Kac--Krein result (see also
De Brange \cite{Bra1:SHS}) on $2\times 2$ canonical systems  with real 
Hamiltonian. Our criteria  cover this as well as some results on 
quasiregularity  from \cite[Sec. 3.2]{KogRof:SIS}. 

Besides, we present several examples which show the limits of the results.

Finally, in Subsection \plref{sec5.5}  we obtain also similar statements on 
quasiregularity of  matrix Sturm-Liouville equation. In the scalar case
these results essentially generalize Krein's result \cite{Kre:TFO} 
(see also \cite{KacKre:SFS}) mentioned above.

In conclusion we mention two recent publications  \cite{Shu:CQC} and          
\cite{Les:ESA} close to our
work (see also references therein) which are devoted to self-adjointness
of elliptic operators on complete manifolds.


\section{The symmetric linear relation induced by a 
first order system}\label{sec2}

In this section we introduce the basic notation about
first order systems. Denote by $\Mat(n,\C)$
the set of complex $n\times n$ matrices and let 
$I\subset\R$
be a (not necessarily open) interval. We denote by $\AC(I)$ 
the set of
all absolute continuous functions on $I, $
i.e. $f\in\AC(I)$ if $f'$ exists a.e., 
is locally integrable, and $f(x)=\int_{x_0}^x f'(s)ds 
+ f(x_0)$.
If $U\subset\R^n$ is an open set, we denote by 
$\AC(I,U)$ the
set of $U$--valued functions whose components lie in 
$\AC(I)$. Finally, if $X$ is a function space over $I$,
then $X_{\comp}$ denotes the subspace consisting of 
those $f\in X$ with compact support in $I$.

With these preparations we consider the first order
system
\begin{equation}\label{G1.2}
J(x)\frac{df}{dx}(x)+B(x) f(x)=\ch(x)g(x), 
     \end{equation}
where $J, B, \ch:I\to \Mat(n,\C)$ are matrix--valued functions
such that:
\begin{align}
&J\in\AC(I,\Mat(n,\C)), &&J(x)=-J(x)^*,\quad \det J(x)\not=0, \;
\textrm{for}\; x\in I,\nonumber\\
       &B\in L^1_{\loc}(I,\Mat(n,\C)), &&
B(x)^*=B(x)-J'(x), \;\textrm{for}\; x\in
       I,\label{G1.0}\\
       &\ch\in L^1_{\loc}(I,\Mat(n,\C)),&& 
\ch(x)=\ch(x)^*, \quad \ch(x)\ge 0,\;\textrm{for}\; x\in
       I.\nonumber
\end{align}
Let $\cLH$ be the set of Borel--measurable  $\C^n$--valued 
functions satisfying
$\scalar{f}{f}_\ch:=\int_I f(x)^*\ch(x)f(x)dx 
<\infty.$
It is well--known (cf. e.g. \cite[Sec. 9]{AkhGla:TLO}, 
\cite{Naj:LDO}) that
$\cLH$ is complete with respect to the semi--norm
$\|f\|_{\ch}=\sqrt{\scalar{f}{f}_\ch}$. Moreover 
$\cLH$ is the completion of $C_{\comp}(I,\C^n)$
with respect to $\|\cdot\|_\ch$. 

We equip $\cLH$ with the (semi--definite) scalar product
   \begin{equation}\label{lr1.1}
\scalar{f}{g}_\ch:=\int_I f(x)^*\ch(x)g(x)dx, 
    \end{equation}
and put
\begin{equation}\LH:=\cLH\Big/\bigsetdef{f\in\cLH}{\|f
\|_\ch=0}.\end{equation}
$\LH$ is a Hilbert space. For a function $f\in \cLH$ 
we will
denote by $\tilde f$ the corresponding class in $\LH$. 
If $\ch(x)$ is invertible a.e. then a class $\tilde f$ 
contains
at most one continuous representative, hence if 
$\ch(x)$ is invertible a.e.
and $f$ is continuous then we 
will not distinguish between $f$ and $\tilde f$.

If in addition $\ch(x)$ is invertible 
for almost all $x\in I$ and $\ch^{-1}, B^*\ch^{-1}B\in 
L^1_{\loc}(I,\Mat(n,\C))$ then \eqref{G1.2} induces a 
symmetric operator
   \begin{equation}
    L:=\ch^{-1}(J \frac{d}{dx}+B)
  \end{equation}
in the Hilbert space $\LH$ with domain $\cd(L)=C^1_{\comp}(I,\C^n)$
(cf. Subsection \plref{sec32} below). 
The symmetry is implied by $B^*=B-J'$ and $\ch^*=\ch$. 
However, the interesting case is the one where $\ch$ is singular.
If $\ch$ is singular then \eqref{G1.2} will in general neither
define an operator nor will it be densely defined. Rather it 
will  give rise to
\emph{symmetric linear relations}, $\csmin$ resp. $S_{\min}$, in  
$\cLH$ resp. $\LH$ as
follows: $\{f, g\}\in \csmin$ if and only if $f\in
\AC_{\comp}(I,\C^n), g\in \cLHc$ and $Jf'+Bf=\ch g$.

For the reader's convenience let us briefly recall the definition of
a symmetric linear relation:
\begin{dfn}\label{SDM-S1.1} Let $\mathfrak H$ be a linear space equipped with a positive
semi--definite hermitian sesqui--linear form $\scalar{\cdot}{\cdot}$.
A linear subspace $\cs\subset \mathfrak{H}\times\mathfrak{H}$ is called
a symmetric linear relation (s.l.r.) if for $\{f_j,g_j\}\in \cs, j=1,2$, one has
$\scalar{f_1}{g_2}=\scalar{f_2}{g_1}$.

For a s.l.r. $\cs$ one defines, as usual, the \emph{domain}
$\cd(\cs):=\setdef{f\in\frH}{\exists_{g\in\frH}\{f,g\}\in\cs}$,
the \emph{range} $\im \cs:=\setdef{g\in\frH}{\exists_{f\in\frH}\{f,g\}\in\cs}$,
and the \emph{kernel} $\ker\cs:=\setdef{f\in\frH}{\{f,0\}\in\cs}$.
Furthermore, the \emph{indeterminant} part of $\cs$ is defined by
$\cs(0):=\setdef{g\in\frH}{\{0,g\}\in\cs}=\ker(\cs^{-1}).$

Finally, the \emph{adjoint} of $\cs$ is $\cs^*:=\setdef{\{f,g\}\in\frH\times\frH}{
\forall_{\{\phi,\psi\}\in\cs} \scalar{f}{\psi}=\scalar{g}{\phi}}$.
\end{dfn}
For example, the graph of an (unbounded) symmetric operator in a Hilbert space
$\mathfrak H$ is a s.l.r. 

$\csmin$ induces a symmetric linear relation, $S_{\min}$, in 
$\LH$ in a fairly
straightforward way: $\{\tilde f,\tilde g\}\in \Smin$ if 
and only if there exist representatives $f\in \tilde f, g\in 
\tilde g$ such that
$\{f,g\}\in \csmin$. Symmetric linear relations arising 
in this way have been studied thoroughly in \cite{Orc:CDE}. 
Unfortunately, \cite{Orc:CDE}
has not been published and therefore is not widely 
available. The
authors received a copy of \cite{Orc:CDE} only after 
the present
work had been almost completed. We emphasize, however, 
that there is
only a small overlap between \cite{Orc:CDE} and the 
present work.

In general $\Smin$ will neither be densely defined nor 
single valued:
\begin{example}\label{ML-S2.2} $I=(0,1), B=0, J=\begin{pmatrix} 0&1\\ 
-1& 0
\end{pmatrix}, \ch(x)=\begin{pmatrix} 1 & 0\\ 
0&0\end{pmatrix}$.
If $\{f,g\}\in \csmin$ then $f_2'=g_1, f_1'=0$, and since
$f$ is continuous with compact support we 
infer $f_1=0$. In view of the special form of $\ch$ this
implies $\tilde f=0$. Hence, the domain of $S$ is $\{0\}$. 
Note that since $g_1=f_2'$ we have $\int_I g_1=0$.

Conversely, given $\tilde g\in\cLH$ with $\int_I g_1=0$
we put $f_2(x):=\int_0^x g_1(s)ds$
and $f_1=0$. Then $\{f,g\}\in\csmin$ and hence
$\{0,\tilde g\}=\{\tilde f,\tilde g\}\in\Smin$.
Consequently, $\Smin=\{0\}\times \setdef{\tilde g}{g\in\cLHc,\;
\int_I g_1=0}$ and $\Smin^*=\setdef{\{\tilde f,\tilde
g\}}{f=\operatorname{const}, g\in\cLH}$.
\end{example}

This example also shows that in general $\Smin$ is 
not closed: 
\begin{dfn}\label{ML-S2.0} We denote by $S$
the closure of $\Smin$, i.e. the minimal closed extension, 
and we put $\Smax:=\Smin^*$. 
Furthermore, we write $\{f,g\}\in\csmax$
if $f,g\in\cLH$, $f$ is absolutely continuous, and 
$Jf'+Bf=\ch g$. Finally, let $\cs$ be the closure of $\csmin$ in $\csmax$,
i.e. $\{f,g\}\in\cs$ if $\{f,g\}\in\csmax$ and
there exists a sequence $(\{f_n,g_n\})_{n\in\N}\subset\csmin$
such that $\|f-f_n\|_\ch, \|g-g_n\|_\ch\to 0, $ as $n\to\infty$.
That is $\cs=\setdef{\{f,g\}\in\csmax}{\{\tilde f,\tilde g\}\in S}$.
\end{dfn}
If $\ch(x)$ is invertible a.e. then $S$ will 
at least be a single valued symmetric operator, i.e. 
$\{\tilde f,\tilde g_1\},
\{\tilde f,\tilde g_2\}\in S$ implies $\tilde 
g_1=\tilde g_2$.
We emphasize that $S$ may be a densely defined 
operator even if $\ch$ is 
singular on a subset of positive Lebesgue measure. 
E. g.  this is the case for $I={\Bbb R}_+$  if 
$\int_\alpha^\beta\ch(t)dt$ is positive definite for 
all 
$\alpha,\beta\in [0,\infty),\ \alpha<\beta$ (see 
\cite{LanTex:GKM}).

A complete description of the indeterminant 
part $S(0)=\setdef{g}{\{0,g\}\in S}$ for $2\times 2$ 
canonical systems has been obtained in \cite{Kac:LRG}, \cite{Kac:LRGC}.

The relations $\cs, S$ will be addressed as the symmetric
linear relation of the first order system 
\eqref{G1.2}. We will write 
\begin{equation}
\cs(J,B,\ch) \quad \text{(resp. } S(J,B,\ch)\text{)} 
\label{ML-G2.6}
\end{equation}
if we want to emphasize the dependence on $J,B,\ch$.

Next we discuss the regularity problem. In view of Definition
\plref{ML-S2.0} integration by parts shows immediately that 
$\{\tilde f,\tilde g\}\in \Smax$ (resp. $S$) if $\{f,g\}\in \csmax$ (resp. $\cs$). 
Denoting by $\pi:\cLH\to \LH$ the quotient map, this 
means that
\begin{equation}(\pi\oplus \pi)(\csmax)\subset \Smax,\quad
(\pi\oplus\pi)(\Smin)\subset S.
\end{equation}
A priori it is not clear whether equality holds. 
We call this the regularity Theorem.

\begin{theorem}[Regularity Theorem] \label{regthm} Let $\{\tilde f,\tilde
g\}\in\Smax$ (resp. $S$).
Then for each representative $g\in\tilde g$ there exists $f\in\tilde f$
such that $\{f,g\}\in\csmax$ (resp. $\cs$).
\end{theorem}

This theorem follows from Propositions \plref{S1.6} and \plref{S1.8} below.
For definite systems (cf. Def. \plref{S1.10} below)
Theorem \plref{regthm} has been proved by Orcutt \cite[Thm. II.2.6 and 
Thm. IV.2.5]{Orc:CDE}. Another proof for (not necessarily definite) 
$2\times 2$ canonical systems was given by
I.S. Kac \cite{Kac:LRGC} in the deposited elaboration 
of \cite{Kac:LRG}.
We note also that his proof is rather long and can not 
be extended 
to $n\times n$ systems. 

In sum, this important regularity result for first 
order systems
is a kind of folklore theorem but proofs are not very 
available in
the literature. To fill this gap and to make this 
article self--contained 
we present a proof below. We emphasize that our
presentation treats the most general case, i.e. we do 
not assume that
the first order system is definite. This is more 
general than \cite{Orc:CDE},
\cite{Kac:LRG}. Also we hope that our presentation is 
simpler and more
perspicuous.

The system \eqref{G1.2}
can be simplified and put into canonical 
form.
The construction is due to Kogan and Rofe--Beketov
\cite[Sec. 1.3]{KogRof:SIS} (see also 
\cite{GohKre:TAV}). 
Since we will make use of 
it heavily and
to fix some notation, let us briefly recall this 
construction:\label{gaugetransf}

A "gauge 
transformation" $U\in \AC(I,\operatorname{GL}(n,\C))$ 
induces a unitary map 
  \begin{equation}\label{2.8}
\Psi_U:\cLH\to \cl_{\tilde \ch}^2(I), \quad f\mapsto 
U^{-1} f,\quad \tilde \ch:= U^*\ch U,
     \end{equation}
and a simple computation shows that 
\begin{equation}
    \Psi_U \cs(J,B,\ch)\Psi_U^*= \cs(\tilde J,\tilde 
B,\tilde \ch),
\end{equation}
where
\begin{equation}\label{ML-G2.8} \tilde J= U^*JU,\quad
    \tilde B= U^*JU'+U^*BU,\quad
    \tilde \ch= U^* \ch U.
\end{equation}
In a first step one chooses $U\in\AC(I,\Mat(n,\C))$ 
such that
$U^* JU=J(0)$. Thus we are reduced to the case where
$J$ is a constant matrix.

In a second step pick $x_0\in I$ and let 
$\fsys(.,\gl):I\rightarrow \Mat(n,\C)$ be the solution 
of the initial
value problem
\begin{equation}
     J\fsys'(x,\gl)+B(x)\fsys(x,\gl)=\gl \ch(x) 
\fsys(x,\gl),\quad \fsys(x_0,\gl)=I_n.
     \label{G1.5}
\end{equation}
Here, $I_n$ denotes the $n\times n$ unit matrix. The 
existence
of $\fsys$ follows from the fact that $B$ and $\ch$ 
are locally integrable.
For $\fsys(x,0)$ we simply write $\fsys(x)$. If $g\in 
\cLHloc$ then, since $\sqrt{\ch}\in \cl^2_{\loc}(I,\Mat(n,\C))$, 
we have
$\ch g\in \cl^1_{\loc}(I,\C^n)$. Thus, the solution of the 
inhomogeneous
initial value problem
\begin{equation}\label{lr1.6}
     J y'(x,\gl)+B(x)y(x,\gl)=\gl \ch(x) y(x,\gl) + \ch(x) g(x), \qquad 
y(x_0,\gl)=0,
\end{equation} 
exists and is unique. Taking into account the well--known
(and easy to verify) formula
   \begin{equation}
\fsys(x,{\ovl{\gl}})^*J\fsys(x,\gl)=J,\quad \gl\in \C,
\label{ML-G2.13}
     \end{equation}
the variation of constants formula reads
\begin{equation}\label{lr1.7}
  y(x,\gl) = (K_\gl g)(x)=\fsys(x,\gl) \int_{x_0}^x 
J^{-1}\fsys(t,{\ovl{\gl}})^* \ch(t)g(t)dt.
\end{equation}
As with $\fsys$ we write $K$ instead of $K_0$.
Now we can choose $Y$ as the gauge transformation. In 
view of 
\eqref{ML-G2.13} and \eqref{ML-G2.8} the gauge 
transformation $Y$
transforms the system into a system 
$\tilde S$ with
\begin{equation}
\tilde J=J(0),\quad \tilde B=0,\quad \tilde \ch=U^*\ch U.
 \label{ML1-G2.15}
\end{equation}
Such systems are called "canonical" in the 
literature.

Another choice of gauge is possible if $\ch$ is 
absolutely continuous
and invertible. Then the gauge $U=\ch^{-1/2}$ turns 
the system into one with $\tilde\ch=1$. 
The interesting cases, however, are those
with singular $\ch$.

Despite the existence of canonical forms obtained from appropriate gauges
we prefer to work in the framework
of \eqref{G1.2} since finding the canonical system corresponding
to the first order system \eqref{G1.2} depends on 
finding the fundamental system of solutions.
Another reason for working in our framework is the 
following:
we will give criteria for $S$ being essentially self--adjoint below.
These criteria are only sufficient and
not gauge invariant, hence it is desirable
to have them at hand also for first order systems 
which are not in canonical
form. It would be nice, however, to have a necessary 
and sufficient
characterization of essential self--adjointness. Such 
a criterion
would necessarily have to be gauge invariant. The 
discovery of such a criterion, however, remains an open problem.

Some remarks are in order about why first order systems are 
interesting. First order systems are not as special
as they seem to be. Namely, an arbitrary symmetric $n^{th}$--order
system is \emph{unitarily} equivalent to a symmetric
first order system (\cite{KogRof:SIS}, \cite{Orc:CDE}). 
In most cases, however, the Hamiltonian $\ch$ of this first order system 
will be singular. Instead of reproducing this result we will
present two important examples. First, we show how a second order 
Sturm--Liouville type (quasi--differential) equation
can be transformed into a system of the form \eqref{G1.2}. 

     \begin{example}\label{SDM-S1.2}
1. We consider a weighted Sturm--Liouville type (quasi--differential) equation
   \begin{equation}
   -\frac{d}{dx}\Bigl(A(x)^{-1}\frac{du}{dx}(x)+Q(x)u(x)\Bigr)+Q(x)^*\frac{du}{dx}(x)+R(x)u(x)=\ch(x) v(x),
\label{SDM-G1.5}
    \end{equation}
where $A,Q,R,\ch\in L^1_{\loc}(I,\Mat(n,\C))$, $A(x)$ is positive definite
for all $x\in I$, and $\ch(x)\ge 0$. The system \eqref{SDM-G1.5} defines a symmetric linear
relation as follows: $\{u,v\}\in\csmin$ if and only if $u\in\AC_{\comp}(I,\C^n),
A^{-1}\frac{du}{dx}+Qu\in \AC_{\comp}(I,\C^n),$
$v\in \cLHc$ and \eqref{SDM-G1.5} holds. 
"Quasi--differential" means that $\frac{du}{dx}$ is not necessarily absolute
continuous.
As for first order systems, let
$\Smin:=\setdef{\{\tilde u,\tilde v\}}{\{u,v\}\in\csmin}$.


Next we introduce the first order system 
 \begin{align}
&\tilde J
\begin{pmatrix} f_1\\f_2\end{pmatrix}'
+\tilde B
\begin{pmatrix} f_1\\f_2\end{pmatrix}
=\tilde\ch
\begin{pmatrix} g_1\\g_2\end{pmatrix},
\label{SDM-G1.6}\\[1em]
&
 \tilde J:=\begin{pmatrix} 0 & iI_n\\iI_n&0\end{pmatrix},\quad \tilde
 B:=\begin{pmatrix} R-Q^*AQ&-iQ^*A\\iAQ&-A\end{pmatrix},
\quad \tilde \ch:=\begin{pmatrix} \ch&0\\0&0\end{pmatrix},
  \label{ML-G2.10}
 \end{align}
and we denote by $\tilde\csmin, \tilde S$ the corresponding s.l.r. in
$\cl_{\tilde\ch}^2(I), L^2_{\tilde\ch}(I)$.

If $\{u,v\}\in\csmin$ then 
$\{(u,i(A^{-1}u'+Qu),(v,0)\}\in\tilde\csmin$.
Conversely, if $\{(f_1,f_2),(g_1,g_2)\}\in\tilde\csmin$ then
$\{f_1,g_1\}\in\csmin$. Hence the unitary isomorphism
\begin{equation}
   \Phi:\LH\longrightarrow L^2_{\tilde\ch}(I),\quad \tilde f\mapsto \tilde{(f,0)}
\end{equation}
implements a unitary equivalence between $\Smin$ and $\tilde \Smin$,
i.e. $(\Phi\times\Phi)^*\tilde\Smin(\Phi\times\Phi)=\Smin$.

Even if
$\Smin$ is (the graph of) a densely defined symmetric operator in the 
Hilbert space $\LH$ the Hamiltonian $\tilde\ch(x)$ is singular everywhere.

2. Consider a general first order system $S=S(J,B,\ch)$ as in \eqref{G1.2}.
We define the square 
of $\csmin$ resp. $\Smin$ as follows:
      \begin{equation}
\begin{split}
     \csmin^2&:=\bigsetdef{\{f,g\}\in\cLH\times\cLH}{\exists_{h\in\cLH}\{f,h\},\{h,g\}\in\csmin},\\
     \Smin^2&:=\bigsetdef{\{\tilde f,\tilde g\}\in\LH\times\LH}{\exists_{\tilde
     h\in\LH}\{\tilde f,\tilde h\},\{\tilde h,\tilde g\}\in \Smin}.
\end{split}
\end{equation}
The squares of $\cs, S$ are defined analogously. We remark first that indeed
\begin{equation}\begin{split}
     \Smin^2&=\bigsetdef{\{\tilde f,\tilde g\}}{\{f,g\}\in\csmin^2},\\
     S^2&=\bigsetdef{\{\tilde f,\tilde g\}}{\{f,g\}\in\cs^2}.
		\end{split}
\end{equation}
To see this consider $\{\tilde f,\tilde g\}\in S^2$ (resp. $\Smin^2$). By definition
there exists a $\tilde h\in\LH$ such that $\{\tilde f,\tilde h\},\{\tilde
h,\tilde g\}\in S$ 
(resp. $\Smin$). Let $g\in\tilde g$. By the regularity Theorem
\plref{regthm} there exists $h\in\tilde h$ such that $\{h,g\}\in \cs$
(resp. $\csmin$, in this case the regularity Theorem is not needed).
Again by the regularity Theorem there exists $f\in\tilde f$ such that
$\{f,h\}\in S$ (resp. $\Smin$). Thus $\{f,g\}\in S^2$ (resp. $\Smin^2$). 
Conversely, if $\{f,g\}\in \cs^2$ (resp. $\csmin^2$) 
then it is clear that $\{\tilde f,\tilde g\}\in S^2$ (resp. $\cs^2$).

Next let  $\{f,g\}\in\csmin^2$, that is there is a $h\in\cLH$ such that 
$\{f,h\}\in\csmin$ and $\{h,g\}\in\csmin$. This is equivalent to the equation
\begin{equation}
\begin{pmatrix} 0 & J\\ J & 0\end{pmatrix}
\begin{pmatrix} f \\ h \end{pmatrix}'
+\begin{pmatrix} 0 & B \\ B & -\ch \end{pmatrix}
\begin{pmatrix} f\\ h\end{pmatrix}
= 
\begin{pmatrix} \ch & 0\\0 & 0\end{pmatrix}
\begin{pmatrix} g \\ 0\end{pmatrix}
\label{ML-G3.2.1}
   \end{equation}
with $f,h\in\AC_{\comp}(I,\C^n), g\in\cLHc$.
A similar argument as under 1. shows that $\Smin^2$ is unitarily equivalent
to $\Smin(J_1,B_1,\ch_1)$, where
\begin{equation}
J_1=\begin{pmatrix} 0 & J\\ J & 0\end{pmatrix},\quad
B_1=\begin{pmatrix} 0 & B \\ B & -\ch \end{pmatrix},\quad
\ch_1=\begin{pmatrix} \ch & 0\\0 & 0\end{pmatrix}.
\end{equation}
Actually, this system is unitarily equivalent
to a system of the form \eqref{ML-G2.10}. Namely, the gauge transformation
\begin{equation}
   U:=\begin{pmatrix} I_n &0\\ 0 & i J^{-1}\end{pmatrix}
   \label{ML-G2.24}
\end{equation}
transforms the system $\Smin(J_1,B_1,\ch_1)$ into $\Smin(\tilde J,\tilde
B,\tilde \ch)$, where
\begin{equation}
\tilde J=\begin{pmatrix}0&i I_n\\i I_n&0\end{pmatrix},\quad
\tilde B=\begin{pmatrix}0&i B^*J^{-1}\\ i J^{-1}B & - (J^{-1})^*\ch
J^{-1}\end{pmatrix},\quad
\tilde\ch=\begin{pmatrix} \ch & 0\\0&0\end{pmatrix}.
\label{ML-G2.25}
\end{equation}
This can be checked using the formulas \eqref{ML-G2.8}.

Note that \eqref{ML-G2.25} is a special case of the structure
\eqref{ML-G2.10}, except that the lower right corner of $\tilde B$
is only positive semi--definite. This is not a surprise since
heuristically $\csmin^2$ can be viewed as a second order system.
\end{example}

For future reference and to fix some notation 
let us present a type of first order systems which
contains the two preceding examples as special cases. Consider
the system
   \begin{equation}
J_1 f'+B_1 f=\tilde \ch g,
\label{ML-G3.2.2}
  \end{equation}
where
  \begin{equation}
J_1=\begin{pmatrix} 0 & J^*\\ -J & 0\end{pmatrix}
,\quad
B_1=\begin{pmatrix} V & B \\ B^*-J' & -A \end{pmatrix},
\quad
\tilde\ch=\begin{pmatrix} \ch & 0\\0 & 0\end{pmatrix}.
\label{ML-G2.27}
  \end{equation}
We assume that \eqref{ML-G3.2.2} satisfies \eqref{G1.0}, that is
$J\in\AC(I,\Mat(n,\C)), V,B,A,\ch\in L^1_{\loc}(I,\Mat(n,\C)$,
$\det J(x)\not=0$, for $x\in I$, $V=V^*, A=A^*$, and $\ch(x)\ge 0$
for $x\in I$. 

As in the previous example, the system \eqref{ML-G3.2.2} can be transformed
quite explicitly onto a system $S(J_2,B_2,\tilde \ch)$ with $J_2$ constant.
We present two normal forms. The gauge transformation \eqref{ML-G2.24}
transforms the system $S(J_1,B_1,\tilde\ch)$ onto $S(J_2,B_2,\tilde\ch)$, where
\begin{equation}
   J_2=\begin{pmatrix} 0&iI_n\\iI_n&0\end{pmatrix},\quad
   B_2=\begin{pmatrix} V&i (B-(J^*)')(J^{-1})^*\\ -i J^{-1} (B^*-J')& - J^{-1}A
   (J^{-1})^*\end{pmatrix}.
\label{ML-G2.26}
\end{equation}
The gauge transformation
\begin{equation}
U:=\begin{pmatrix} I_n& 0\\ 0& i I_n\end{pmatrix}
\end{equation}
transforms the system $S(J_2,B_2,\tilde\ch)$ onto $S(J_3,B_3,\tilde \ch)$,
where
\begin{equation}
    J_3=\begin{pmatrix} 0&-I_n\\I_n&0\end{pmatrix},\quad
    B_3=\begin{pmatrix} V&-(B-(J^*)')(J^{-1})^*\\ -J^{-1} (B^*-J')& - J^{-1}A
   (J^{-1})^*\end{pmatrix},
\label{ML-G2.26a}
\end{equation}
Note that the normal form \eqref{ML-G2.26} as well as \eqref{ML-G2.26a} are
special cases of \eqref{ML-G2.27}.

These systems will serve as a source of examples and they will be discussed
at several places through the course of the paper.

\subsection{The finite interval case, regularity}\label{sec21}

In this subsection we consider a finite interval 
$I=(a,b)$,
$-\infty<a<b<\infty$. Moreover, we assume that $\ch, 
B\in \cl^1(a,b)$. In view of the previous discussion of
gauge transformations w.l.o.g. we may assume that
$J(x)=J(0)=:J$ is constant.
We denote by $\fsys(.,\gl)$ the solution of 
\eqref{G1.5} with
$x_0=a$.

We introduce the linear map
\begin{equation}\begin{split}
    \delta_\gl:&\cLH\longrightarrow \C^n,\\
         &g\mapsto J\fsys(b,\gl)^{-1}(K_\gl g)(b)
                          =\int_a^b 
\fsys(t,{\ovl{\gl}})^*\ch(t) g(t) dt.
		\end{split}
\end{equation}
Obviously, $\delta_\gl$ induces a map on $\LH$. We 
will be sloppy here and do not distinguish 
between $\delta_\gl$ and its 
induced map on $\LH$.
For $\delta_0$ we just write $\delta$.
Note that since $\delta_\gl$ is continuous and since 
the target space
$\C^n$ is finite--dimensional we have
\begin{equation}
    \im \delta_\gl= 
\delta_\gl(\cl^2_{\ch,\comp}(I)).
\end{equation}
We have even more:
\begin{lemma} \label{S1.3}
$\cLHc\cap\ker\delta_\gl$ is dense in $\ker 
\delta_\gl$.
\end{lemma}
\begin{proof} Let $g_1,\ldots, g_k\in\cLHc$ such that
$\delta_\gl(g_1),\ldots,\delta_\gl(g_k)$ is a basis of 
$\im\delta_\gl$.
Then we have topological direct sum splittings
\begin{equation}\begin{split}
    \cLHc&=(\cLHc\cap\ker\delta_\gl)\oplus 
<g_1,\ldots,g_k>,\\
    \cLH&=\ker\delta_\gl\oplus <g_1,\ldots,g_k>.
		\end{split}
\label{G1.11}
\end{equation}
This implies the claim.\end{proof}

\begin{cor}\label{S1.3a} Let $\{f,g\}\in\csmax$. Then, 
for
$\{\tilde f,\tilde g\}$ to be in $S$ it is 
sufficient
that $f(a)=f(b)=0$. 
\end{cor}
\begin{proof}
$f(a)=f(b)=0$ implies $g\in\ker\delta$ and, 
in view of the previous lemma, we may choose a 
sequence
$(g_n)\subset\ker \delta\cap\cLHc$
with $g_n\to g$ in $\cLH$. Then $Kg_n\in 
\AC_{\comp}(I,\C^n)$
and $Kg_n\to Kg=f$ in $\cLH$. Thus 
$\{\widetilde{Kg_n},\tilde g_n\}\in \Smin$
and $\{\widetilde{Kg_n},\tilde g_n\}\to \{\tilde 
f,\tilde g\}$.
\end{proof}

We put 
\begin{equation}
   \Mmat(\gl)=\Mmat(\cs,\gl)=\int_a^b 
\fsys(x,\gl)^*\ch(x)\fsys(x,\gl)dx.
   \label{G1.12}
\end{equation}
For $\Mmat(0)$ we just write $\Mmat$.

\begin{lemma}[cf. {\cite[Thm. 
1.1]{KogRof:SIS}}]\label{S1.4} 
$\ker \Mmat(\gl), \im \Mmat(\gl)$ are independent 
of $\gl$, in particular $\rank \Mmat(\gl)$ is 
independent of $\gl$.
    \end{lemma}
\begin{proof} Fix $\gl_0,\gl\in\C$ and consider 
$\xi\in\ker \Mmat(\gl)$.
Then we have
\begin{equation}
    \int_a^b \xi^*\fsys(x,\gl)^*\ch(x)\fsys(x,\gl)\xi 
dx=0
\end{equation}
and hence $\ch(x)\fsys(x,\gl)\xi=0$ for almost all 
$x\in I$.
Moreover, the function $f(x)=\fsys(x,\gl)\xi$ 
satisfies the differential equation
\begin{equation}
      J f'(x)+Bf(x)=\gl \ch(x) f(x)=\gl_0\ch(x)f(x)
\end{equation}
for almost all $x\in \R$. Thus, by the uniqueness 
theorem for first
order differential equations we have 
$f(x)=\fsys(x,\gl_0)f(a)=\fsys(x,\gl_0)\xi$.
Moreover, since $\xi\in\ker \Mmat(\gl)$,
\begin{equation}
      0=\xi^*\Mmat(\gl)\xi=\int_a^b f(x)^*\ch(x) 
f(x)dx=\xi^* \Mmat(\gl_0)\xi.
\end{equation}
Since $\Mmat(\gl_0)\ge 0$ we infer $\xi\in\ker 
\Mmat(\gl_0)$. 

Since $\gl_0,\gl$ were arbitrary we have proved that 
$\ker \Mmat(\gl)$ is
independent
of $\gl$. This implies the rest of the assertions.
\end{proof}

The rank of $\Mmat$ will play a crucial role, thus we 
put
\begin{equation}
    \rank(\cs):=\rank(S):=\rank(\Mmat).
    \label{G1.16}
\end{equation}

\begin{lemma} \label{S1.5} 
$\im \delta_\gl=\im 
\Mmat=\setdef{\xi\in\C^n}{\ch \fsys\xi=0
 \,\textrm{a.e.}}^\perp$.

Moreover, we have an orthogonal sum decomposition
    \begin{equation}
        \cLH=\ker\delta_\gl\oplus 
\bigsetdef{\fsys(.,\ovl{\gl})\xi}{\xi\in\im \Mmat}.
        \label{G1.17}
   \end{equation}
\end{lemma}
    \begin{proof} First we prove \eqref{G1.17}. For any
$\xi\in\C^n$ and $g\in\cLHc$ one has
\begin{equation}\begin{split}
 \scalar{\xi}{\delta_\gl(g)}&=\int_a^b 
\xi^*\fsys(x,\ovl{\gl})^*\ch(x)g(x)dx\\
           &=\int_a^b\big(\fsys(x,\ovl{\gl})\xi)^* 
{\ch}(x)g(x) dx =
           \scalar{\fsys(.,\ovl{\gl})\xi}{g}_{\ch},
		\end{split}
     \end{equation}
hence $\delta_\gl^*(\xi)=\fsys(.,\ovl{\gl}) \xi.$ We 
note that  $\fsys(.,\ovl{\gl})\xi=0$ in $\LH$  
(that is $\ch\fsys(.,\ovl{\gl})\xi=0)$ for $\xi \in 
(\im \delta_\gl)^\perp$. Thus one infers
    \begin{equation}\label{G1.17a}
\cLH=\ker\delta_\gl\oplus\im \delta_\gl^* = 
          \ker\delta_\gl\oplus 
\bigsetdef{\fsys(.,\ovl{\gl})\xi}{\xi\in\im 
\delta_\gl}.
       \end{equation}
It follows that each $g\in\cLHc$ admits a unique 
decomposition
   \begin{equation}    
  g=g_0+Y(\cdot,\ovl{\gl})\xi_g, \qquad 
g_0\in\ker\delta_\gl,\ \ 
\xi_g\in\im\delta_\gl,
  \label{G1.19}
       \end{equation}
where $\xi_g$ is the unique element
in $\im \delta_\gl$ such that 
$\delta_\gl(\fsys(.,\ovl{\gl})\xi_g)=\delta_\gl(g)$.
Furthermore, 
    \begin{equation*}
\delta_\gl(g) = \int_a^b 
\fsys(x,\ovl{\gl})^*\ch(x)\fsys(x,\ovl{\gl}))\xi_gdx =
    \Mmat(\ovl{\gl})\xi_g,\qquad g\in \cLH.
    \end{equation*}
Hence $\im\delta_\gl \subset \im 
\Mmat(\ovl{\gl})=\im \Mmat.$ 
Since the opposite inclusion is 
obvious one gets $\im\delta_\gl = \im \Mmat.$ In 
view of (\ref{G1.17a}) this
relation implies (\ref{G1.17}). To complete the proof 
it remains to note that
$\ker \Mmat = \setdef{\xi\in\C^n}{\ch \fsys\xi=0 
\,\textrm{a.e.}}.$
        \end{proof}
  \begin{prop}\label{S1.6}
\begin{thmenum}
\item For all $\gl\in\C$ we have
\begin{align*}
           \im (\Smax-\lambda)&= \LH,\\
           \im (S-\lambda)&=\pi (\ker\delta_\gl)=\bigsetdef{\pi g}{g\in
           \cLH, 
\int_a^b\fsys(x,\ovl{\gl})^*\ch(x)g(x)dx=0},\\
           \ker(S-\lambda)&=\{0\},\\
           \ker(\Smax-\lambda)&=\bigsetdef{\pi 
\fsys(.,\lambda)\xi}{\xi\in \im \Mmat}
           \simeq \im \Mmat.
      \end{align*}
\item If $\{\tilde f,\tilde g\}\in \Smax$ then for 
each representative
           $g\in\tilde g$ there exists $f\in \tilde 
f$, $f\in \AC(I,\C^n)$,
such that $Jf'+Bf=\ch g$. 
In particular 
$\pi_2(\csmax):=(\pi\oplus\pi)(\csmax)=\Smax$.
\item $\pi_2\bigl(\setdef{\{f,g\}\in \csmax}{f(a)=f(b)=0}\bigr)=S$. 
Moreover, 
   \begin{equation*}
\cs= (\pi_2^{-1}S)\cap {\cs}^* = 
\bigsetdef{\{f,g\} \in{\cs}^*}{\ f(a)\in \ker \Mmat,\  
f(b)=\fsys(b)f(a)}.
     \end{equation*}
\end{thmenum}
\end{prop}
\begin{proof} 
(1) If $g\in \cLHc$ is
arbitrary then $\{\widetilde{K_\gl g},\tilde g\}\in 
(\Smax-\gl)$ and we have
proved that $\im (\Smax-\gl)=\LH$.

If $g\in\ker\delta_\gl$ then by Corollary 
\plref{S1.3a} we have
$\{\widetilde{K_\gl g},\tilde g\}\in (S-\gl)$, 
thus 
$\ker \delta_\gl\subset \im(S-\gl).$
Since $\im (\Smin-\gl)\subset\ker\delta_\gl$
by definition and since $\delta_\gl$ is continuous we 
conclude that
$\im (S-\gl) \subset \ovl{\im (\Smin-\gl)} 
\subset\ker\delta_\gl$. 
We have proved $\ker\delta_\gl=\im(S-\gl)$.
Furthermore we infer $\ker(S-\gl)=\im (\Smax-
\ovl{\gl})^\perp=\{0\}$ and
$\ker (\Smax-\gl)=\im (S-
\ovl{\gl})^\perp=(\ker\delta_{\ovl{\gl}})^\perp
=\setdef{\pi\fsys(.,\gl)\xi}{\xi
\in\im \Mmat},$ in view of \eqref{G1.17}. 

(2) Let $\{\tilde f,\tilde g\}\in \Smax$ and let 
$f\in\tilde f,g\in\tilde g$.
We put $f_1(x):=Kg(x)$. Then $\{\tilde f-\tilde 
f_1,0\}\in \Smax$, i.e.
$\tilde f-\tilde f_1\in \ker \Smax$. Consequently, 
there is a $\xi\in\im \Mmat$
such that $\tilde f=\tilde f_1+\widetilde{\fsys\xi}$ 
and hence
$f_2:=f_1+\fsys\xi$ is an absolute continuous 
representative of $\tilde f$
which satisfies $Jf_2'+Bf_2=\ch g$.

(3) Let $\{\tilde f,\tilde g\}\in \Smax$ with 
representatives
$\{f,g\}\in\csmax$. Then
\begin{equation}
      f(x)=\fsys(x) f(a)+Kg(x).
\end{equation}
If $\{\tilde f,\tilde g\}\in S$ then by (1) we have
$g\in \ker\delta$ and
hence $f(b)=\fsys(b)f(a)$. Moreover, 
$\{\widetilde{Kg},\tilde g\}\in S$
and thus $\{\widetilde{\fsys f(a)},0\}\in\ker 
S=\{0\}$.
This implies $\ch\fsys f(a)=0$ a.e. and thus 
$f(a)\in\ker \Mmat$.

Conversely, let $f(b)=\fsys(b)f(a)$ and $f(a)\in\ker 
\Mmat$. Then $Kg=f-\fsys f(a)$
represents the same element $\tilde f\in \LH$ as $f$. 
Moreover
$f(b)=\fsys(b)f(a)$
implies $\delta(g)=0$, hence $\{\tilde f,\tilde 
g\}=\{\widetilde{Kg},\tilde
g\}$. Since $Kg(a)=Kg(b)=0$ this argument also shows 
$\pi_2\bigl(\setdef{\{f,g\}\in \csmax}{f(a)=f(b)=0}\bigr)=S$.
\end{proof}

\subsection{Arbitrary intervals}\label{sec22}

Now we consider an arbitrary, finite or infinite, 
interval $I\subset\R$. Let
$J,B,\ch$ be as in \eqref{G1.0} with $J=J(0)$ constant. We fix a point
$x_0\in I$ and denote by $\fsys(x,\gl)$ the solution 
\eqref{G1.5}. For any finite
subinterval $\tilde I\subset \Icirc$, $I^\circ:=I\setminus \partial I$, we consider the matrix
\begin{equation}
    \Mmat_{\tilde I}(\gl):=\int_{\tilde I} 
\fsys(x,\gl)^*\ch(x)\fsys(x,\gl)dx.
    \label{G1.21}
\end{equation}
In view of Lemma \plref{S1.4} the range of 
$\Mmat_{\tilde I}(\gl)$ is independent
of $\gl$ and as before we write $\Mmat_{\tilde I}$ 
instead of $\Mmat_{\tilde I}(0)$.
Note, however, that $\Mmat_{\tilde I}(\gl)$ depends on 
the choice of the base point $x_0$.
$\tilde I\mapsto \Mmat_{\tilde I}$ is an increasing 
map with values in
the positive semi--definite matrices. Moreover, in view 
of \eqref{G1.21}
$\Mmat_{\tilde I}$ depends continuously on the 
endpoints of $\tilde I$. Since
the rank is a lower semi--continuous function on the 
space of $n\times n$
matrices we infer that there exists a compact interval 
$I_0\subset \Icirc$
such that for any compact interval $I_0\subset \tilde 
I\subset \Icirc$
we have 
\begin{equation}
\im \Mmat_{I_0}=\im \Mmat_{\tilde 
I}.\label{G1.22}
\end{equation}
We then put (cp. \eqref{G1.16})
\begin{equation}
      \rank(S):=\rank(\cs):=\rank \Mmat_{I_0}.
\label{G1.23}
\end{equation}
Somewhat sloppy, in view of \eqref{G1.22}, we will 
write $\ker \Mmat, \im \Mmat$ for $\ker 
\Mmat_{I_0}, \im
\Mmat_{I_0}$. For $g\in \cl^2_{\ch,\comp}(\Icirc)$ we put
\begin{equation}
\delta_\gl(g):=\int_I\fsys(x,\ovl{\gl})^*\ch(x)g(x) 
dx.
  \label{G1.24}
\end{equation}

\begin{lemma}\label{S1.7} Let $k=\rank(S)$. Then there 
exist
$g_1,\ldots,g_k\in\cLHc$ such that 
there is a direct sum decomposition
\begin{equation}
     \cLHc=\ker\delta_\gl\oplus <g_1,\ldots,g_k>.
\end{equation}
\end{lemma}
\begin{proof}  In view of Lemma \plref{S1.5} and the 
previous considerations
we have $\im \delta_\gl=\im \Mmat_{I_0}$. Hence, 
from \eqref{G1.11} we infer
that we may choose $g_1,\ldots,g_k\in 
\cl^2_{\ch,\comp}(I_0)$ such that
$\delta_\gl(g_1),\ldots,\delta_\gl(g_k)$ is a basis of 
$\im \delta_\gl$. 
This implies the assertion.
\end{proof}

Now we are in the position to prove the analogue of 
Proposition \plref{S1.6} for general intervals.

\begin{prop}\label{S1.8} 
Let $\cs$ be the symmetric linear relation induced by
the first order system \eqref{G1.2} on an arbitrary 
interval
$I$. Then:
\begin{thmenum}

\item  $\im(\csmin-\gl)\supset \ker \delta_\gl$. 
Moreover, if $I=[0,b)$ is left--closed (resp. $I=(a,0]$ right--closed)
then $\im (\Smax-\gl)$ is dense in $\LH$ and 
$\ker (S-\gl)=\{0\}$.
\item If $\{\tilde f,\tilde g\}\in \Smax$ then for 
each representative
           $g\in\tilde g$ there exists $f\in \tilde 
f$, $f\in \AC(I,\C^n)$,
such that $Jf'+Bf=\ch g$. In particular 
$\pi_2(\csmax)=\Smax$.
\item  Let $I=\R_+$ and let ${\Mmat}_0$ and ${\Mmat}_1$ be
the matrices constructed in
\eqref{G1.21}--\eqref{G1.23}
with respect to the base point $c\in [0, \infty]$ and the 
intervals $[0,c]$ and $[c,\infty)$ respectively.
Suppose also that $\im ({\Mmat}_0)=\im ({\Mmat}_1).$
Then for each $\xi \in \im (J^{-1}\Mmat)$
there exists $\{f,g\}\in \csmax$ with compact support 
such that $f(c)=\xi$. Moreover, $\{f,g\}\in \csmin$ if $c>0.$  
\item Let $I=\R_{\pm}$ and 
let $\{\tilde f,\tilde g\}\in S$ with 
representatives
$\{f,g\}\in\csmax$. Then $f(0)\in \ker \Mmat$. 
Moreover, $\pi_2\bigl(\setdef{\{f,g\}\in\cs}{f(0)=0}\bigr)=S$.
\end{thmenum}
   \end{prop}
\begin{proof} For simplicity we will give the proof 
for $\gl=0$.

(1) Let $I=[0,b)$ be left--closed and let $g\in \cLHc$. Then
choose $c>\max(\supp g)$ and put
\begin{equation}
 f(x):=\fsys(x)\int_a^x J^{-1}\fsys(t)^* \ch(t) 
g(t)dt.
 \label{G1.25}
\end{equation}
Since $I$ is left--closed we then have $\{\tilde f,\tilde g\}\in \Smax$ and hence 
$\LHc\subset\im \Smax$. Thus
$\Smax$ has dense range and consequently $\ker S=\{0\}$.

The same construction shows for any interval $I$ 
that if $g\in\ker\delta\cap\cl^2_{\ch,\comp}(\Icirc)$ 
then the function $f$ has compact support in $\Icirc$ and thus 
$\im\csmin\supset\ker\delta$.

(2) Let $f_0\in\tilde f$ be any representative and 
put $f_1(x)=\fsys(x)\int_{x_0}^xJ^{-
1}\fsys(y)^*\ch(y)g(y)dy$. Then
$f_1$ is absolutely continuous. Using integration by 
parts and
\eqref{ML-G2.13} one obtains for
any pair $\{\varphi,\psi\}\in \csmin$
\begin{equation}
      \int_I f_1^*\ch \psi=\int_I g^*\ch 
\varphi=\int_I f_0^*\ch\psi. 
\label{G1.26}
\end{equation}
By (1) we have $\im\csmin\supset\ker\delta$, thus \eqref{G1.26} implies
\begin{equation}
        \int_I (f_0-f_1)^*\ch \psi=0,\quad \text{for all}\quad \psi\in\ker\delta.
\label{G1.27}
\end{equation}
Since the $g_j$ in Lemma 
\plref{S1.7} satisfy
$\supp(g_j)\subset I_0$ we apply
Lemma \plref{S1.5} and Lemma \plref{S1.7} to conclude 
that there is a $\xi\in\im \Mmat$ such that
for all $\psi\in\cLHc$ one has
\begin{equation}
        \int_I(f_0-f_1-\fsys\xi)^*\ch\psi=0.
\label{G1.28}
\end{equation}
Note that by integration by parts one has 
$\int_I\xi^*\fsys(x)^*\ch(x)u(x)dx=0$
for all $u\in\ker\delta$, even if $\supp(u)\cap 
(I\setminus I_0)\not=\emptyset$.
\eqref{G1.28} implies that $f=f_1+\fsys\xi$ is an 
absolute continuous
representative
of $\tilde f$ with $J f'+Bf=\ch g$.

(3) We may assume that $I_0=[a_0, c]$ and $I_1=[c, a_1]$ where $a_0>0$. 
Then choose $\eta_0, \eta_1  \in {\mathbb C}^n$ satisfying 
$\xi =J^{-1}\Mmat_0\eta_0=-J^{-1}\Mmat_1\eta_1$ and put
    \begin{equation*}
g(t)=\begin{cases}
\chi _0(t) \fsys(t) \eta _0, &  t\in [0,c),\\
\chi _1(t) \fsys(t) \eta _1, &  t\in [c,\infty ).
     \end{cases}
      \end{equation*}
Here $\chi_0$ and $\chi_1$ are  the 
characteristic functions of the intervals
$I_0=[a_0, c)$ and $I_1=[c, a_1]$ respectively. 
Then we define $f$ by (\ref{G1.25}) with $a$ replaced by $a_1.$
It is clear that $\supp f\subset [0, a_1]$ and 
      \begin{equation*}
f(c)=J^{-1}\int^c_{a_1}\fsys 
(t)^*{\ch}(t)g(t)dt=
-J^{-1}\Mmat_1\eta_1 =\xi.
    \end{equation*}
Furthermore, for $x\in [0, a_0]$ one gets 
      \begin{align*}
f(x)&=\fsys (x)J^{-1} \int ^x_{a_1}\fsys^*(t){\ch}(t)g(t)dt\\
&=-\fsys (x)J^{-1}\Bigl[\int ^c_{a_0}\fsys^*(t){\ch}(t)\fsys(t)dt\,\eta _0 
+\int ^{a_1}_c\fsys^*(t){\ch}(t)\fsys(t)dt\,\eta _1\Bigr]\\
&=-\fsys(x)J^{-1}[\Phi_0\eta_0 + \Phi_1\eta_1]=0.
     \end{align*}

(4) Let $\xi \in  \im J^{-1}\Mmat$. According to (3)
we may choose $\{\varphi, \psi \}\in \csmax$ with 
compact support 
such that $\varphi (0)=\xi$. For each 
$\{\tilde f,\tilde g\}\in  S$ we have on the one 
hand 
$(\varphi, g)_{{\ch}}=(\psi, f)_{{\ch}}.$ 
Since $\varphi,\psi$ have compact support we may 
integrate by parts and
thus find
       \begin{equation*}
0=\varphi(0)^*J f(0)=-
\scalar{J\varphi(0)}{f(0)}=\scalar{J\xi}{f(0)}.
        \end{equation*}
Thus $f(0)$ is orthogonal to $\im \Mmat$, that is 
$f(0)\in \ker \Mmat.$

To prove the last assertion let $\{\tilde f,\tilde g\}\in S$ with
representatives $\{f,g\}\in\cs$. Then $f(0)\in\ker\Phi$ and
hence $f_1:=f-Yf(0)$ is an absolute continuous representative of
$f$. Moreover, $Jf_1'+Bf_1=g$ and $f_1(0)=0$. Consequently,
$\{\tilde f,\tilde g\}=\{\tilde f_1,g\}$ and $\{f_1,g\}\in S$.
  \end{proof}
     \begin{remark} The converse of (4) does not hold 
without further assumptions.
Roughly speaking the system has to be 
"in the limit point case" at infinity.
We will give criteria under which this is true. 
  \end{remark}

    \subsection{Definite first order systems}\label{sec23}

   \begin{dfn}\label{S1.10} The system \eqref{G1.2} is 
said to be
\emph{definite} on $I$ 
if $\ker \Mmat=\{0\}$.
In other words there is a compact subinterval 
$I_0\subset \Icirc$ such that
for all intervals $I_0\subset\tilde I\subset \Icirc$ and 
all $\gl\in \C$ the
matrix $\Mmat_{\tilde I}(\gl)$ (cf. \eqref{G1.21}, 
Lemma \plref{S1.4}) is invertible.
   \end{dfn}

In other words, the system \eqref{G1.2} is definite if $0$ is the only solution of
\[ Jf'+Bf=0,\quad \ch f=0\]
in $\cLH$.

The property of a system \eqref{G1.2} to be definite 
is gauge invariant. 
For a canonical system ($J=J(0), B=0$) this 
property may be 
reformulated solely in terms of the Hamiltonian $\ch.$ 
Namely, it is shown in \cite{GohKre:TAV} and 
\cite{KogRof:SIS} 
that a canonical system is definite iff 
the Hamiltonian 
$\ch$ is of positive type, that is $\int_{I_0}\ch$ is 
invertible for some 
$I_0.$ 

Note also that the system \eqref{G1.2} is definite for 
arbitrary 
$B$ and $J$ if the Hamiltonian $\ch$ is positive 
definite on a subset
of positive Lebesgue measure. 
We emphasize however that for a 
general system \eqref{G1.2} being definite is a property of
the system and depends on $J, B,$ too.
The two examples show that 
the invertibility of $\int_{I_0}\ch$ is unrelated to 
$\ch$ being of
positive type.

The usefulness of the notion of definiteness mainly stems from the
following fact:

\begin{prop}\label{ML-S2.15} Assume that the system \eqref{G1.2} is definite.
Let $\{\tilde f,\tilde g\}\in\Smax$. Moreover let $\{f_j,g_j\}\in\csmax, j=1,2$,
be representatives of $\{\tilde f,\tilde g\}$, i.e. $f_j\in\tilde f,
g_j\in\tilde g$.
Then $f_1=f_2$. 
\end{prop}
\begin{remark}\label{ML-S2.15a} Note that the Proposition does not say that $\tilde f$ has
exactly one absolute continuous representative. In fact it is easy to see that
this is false. See the third example below. 

However, Proposition \plref{ML-S2.15} allows to speak of the value of $\tilde
f$ at a point. I.e. for $x\in \R$ put $\tilde f(x):=f(x)$, where
$\{f,g\}\in\csmax$ is a representative of $\{\tilde f,\tilde g\}$.
Proposition \plref{ML-S2.15} says that $\tilde f(x)$ is well--defined
independently of the choice of $\{f,g\}$.
\end{remark}
\begin{proof} Consider $\{\varphi,\psi\}:=\{f_1-f_2,g_1-g_2\}\in\csmax$. Then
$\{\tilde\varphi,\tilde\psi\}=0$ and hence
\begin{equation}
   J\varphi'+B\varphi=0,\quad \ch\varphi=0.
\end{equation}
Then the definiteness implies $\varphi=0$ and we are done.
\end{proof}

\begin{example}\label{Mark-S1.11}
1. Let 
\begin{equation}
J=\begin{pmatrix} 0&1\\-1&0\end{pmatrix},\quad B=-
I_2,\quad
\ch(x)=\begin{pmatrix} \cos^2(x)& \sin(x)\cos(x)\\ 
\sin(x)\cos(x)&
\sin^2(x)   \end{pmatrix},
\end{equation}
and $I=[0,\pi]$. Then 
$\int_0^\pi\ch(x)dx=\frac{\pi}{2}I_2$
is invertible. However, the function
\begin{equation}
             f(x) =\begin{pmatrix}\sin(x)\\ -
\cos(x)\end{pmatrix}
\end{equation}
satisfies $Jf'+Bf=0$ and $\ch f=0$. Thus the system is 
not definite.

Note, that for this system we have
\begin{equation}
    \fsys(x)=\begin{pmatrix} \cos(x) & -\sin(x)\\     
                             \sin(x) & \cos(x)
              \end{pmatrix}.
\end{equation}
Using this as gauge (cf. \eqref{ML1-G2.15}) we 
obtain the corresponding
canonical system $S(\tilde J, \tilde B,\tilde{\ch})$ 
with 
$\tilde J=J, \tilde B=0, \tilde 
\ch=\diag(1,0)$.
It is clear that this system is not definite.

2. Let $V\in L^1(I)$ and put
\begin{equation}
J=\begin{pmatrix}0&-1\\1&0\end{pmatrix},\quad 
B=\begin{pmatrix}
V&0\\0&-1\end{pmatrix},
\quad \ch(x)=\begin{pmatrix} 1&0\\0&0\end{pmatrix}.
\end{equation}
Then it is easy to check that $\LH\simeq L^2(I)$, and 
the equation
$Jf'+Bf=\ch g$ is equivalent to $-f_1''+Vf_1=g_1$. 
This shows that
the system is equivalent to the Schr{\"o}dinger 
operator
$-\frac{d^2}{dx^2}+V$ on the interval $I$.

Now assume that $I$ is a finite interval. Then 
$\int_I\ch$ is
of rank one and hence not invertible. We claim, 
however, that
the system is definite. Namely, let $Jf'+Bf=0$
and $\int_If^*\ch f=0$. Then $f_1=0$ and since 
$f_2=f_1'$ we also
have $f_2=0$.

Another way of seeing this is to look at the 
fundamental system $\fsys$.
$\fsys$ is a Wronski matrix
\begin{equation}\fsys=\begin{pmatrix}
                   f & g\\
                   f'& g'
		\end{pmatrix},
\end{equation}
thus
\begin{equation}\tilde \ch=\fsys^*\ch 
\fsys=\begin{pmatrix}
                                   f^2 & fg\\
                                   fg  & g^2
					    \end{pmatrix}.
\end{equation}
Since $f,g$ are linearly independent the 
Cauchy--Schwarz--Bunyakovskii inequality yields
$\int_I\tilde\ch>0$.
This example is a special case of Example \plref{SDM-S1.2}.
See also Proposition \plref{ML-S2.16} for a more general
result on definiteness.

3. In 2. consider the special case $V=0$. Put $f:=\binom{1}{0}, g=0$. Then
$\{f,g\}\in\csmax$. However, $\binom{1}{1}$ is a second absolute continuous
representative of $\tilde f$. This is an example for the claim made
in Remark \plref{ML-S2.15a}.
\end{example}

The last example is a special case of the following definiteness result
for systems of the form \eqref{ML-G2.26}, \eqref{ML-G2.26a}.
      \begin{prop}\label{ML-S2.16} Let $I\subset\R$ be an interval.
We consider the system $S(\tilde J,\tilde B,\tilde \ch)$, where
\begin{equation}
   \tilde J=\begin{pmatrix} 0&-I_n\\I_n&0\end{pmatrix},\quad
   \tilde B=\begin{pmatrix} V& B\\B^*&-A\end{pmatrix},\quad
   \tilde \ch=\begin{pmatrix} \ch &0\\0&0\end{pmatrix}
\end{equation}
are as in \eqref{ML-G2.26a}. Assume that the set
$I_0:=\setdef{x\in I}{\det(A(x)\ch(x))\not=0}$ has positive
Lebesgue measure. Then the system $S(\tilde J,\tilde B,\tilde \ch)$
is definite.
       \end{prop}
\begin{proof} Consider $f\in\cl^2_{\tilde\ch}(I)\cap\AC(I,\C^{2n})$ satisfying
\begin{equation}
   \tilde J f'+\tilde Bf=0,\qquad \int_I f^*\tilde \ch f=0.
\label{ML-G2.44}
\end{equation}
We have to show that $f=0$. \eqref{ML-G2.44}
translates into
\begin{align}
      f_1'+ B^*f_1-A f_2&=0,\label{ML-G2.45a}\\
      -f_2'+Bf_2+V f_1  &=0,\label{ML-G2.45b}\\
          \int_{I} f_1^*\ch f_1&=0.\label{ML-G2.45c}
\end{align}
\eqref{ML-G2.45c} implies that $\ch f_1=0$ a.e. Thus the set
$I_1:=\setdef{x\in I}{\det(A(x)\ch(x))\not=0, \ch(x)f_1(x)=0}=
\setdef{x\in I}{\det(A(x)\ch(x))\not=0, f_1(x)=0}$ has positive
Lebesgue measure. A set of positive Lebesgue measure contains an accumulation
point
of itself; the reason is that a subset of the reals which does not contain an accumulation
point of itself is at most countable. So let $x_0\in I_1$ be an accumulation
point of $I_1$. Then $f_1(x_0)=f_1'(x_0)=0$ and by \eqref{ML-G2.45a}
$A(x_0)f_2(x_0)=0$. Since $A(x_0)$ is invertible we infer $f_1(x_0)=f_2(x_0)=0$
and hence $f(x_0)=0$. Since $f$ is a solution of the homogeneous first order
equation $\tilde f'+\tilde B f=0$ this implies $f=0$.
\end{proof}

  \subsection{Formal defect subspaces}\label{sec24}
In this section we present some results on the 
square--integrable solutions of
the system 
\begin{equation}\label{G2.43}
J(x)y'(x)+B(x) y(x)=\gl\ch(x)y(x). 
     \end{equation}

Let
\begin{equation}\begin{split}
      \ce_{\gl}(S)&:=\bigsetdef{f\in\cLHv{I}\cap 
\AC(I,\C^n)}{ Jf'+ 
Bf=\gl \ch f}\\
   &=\bigsetdef{f\in\cLH}{\{f,\gl f\}\in \cs^*}\\
   &=\ker(S^*-\gl),
		\end{split}
       \label{G1.33}
\end{equation}
and denote by
\begin{equation}
     \cn_{\pm}(S):=\dim \ce_{\pm i}(S)
\end{equation}
the \emph{formal} deficiency indices of the system 
\eqref{G1.2}. Furthermore, for a symmetric linear relation $A$ in 
the Hilbert space $\mathfrak H$ we denote by
   \begin{equation}\begin{split}
        E_\gl(A)&:=\bigsetdef{f\in \mathfrak 
H}{\{f,\gl f\}\in A^*}\\
     &=\ker(A^*-\gl)
		   \end{split}
,\qquad \gl\in\C,
\label{SDM-G1.13}
   \end{equation}
the defect subspace and by
\begin{equation}
       N_\pm(A):=\dim E_{\pm i}(A)
\end{equation}
the deficiency indices of $A$. 
It is well--known (see \cite{AkhGla:TLO}, 
\cite{Naj:LDO}) that  
\begin{equation}
\dim E_{\pm \gl}(A)=N_{\pm}(A),\qquad \gl \in \C_+:=\bigsetdef{z\in\C}{\Im
z>0}.
\label{SDM-G1.14}
\end{equation}
 We present however two simple proofs of \eqref{SDM-G1.14}.

The first proof follows from the observation that the 
relation $A^*-\lambda$
is semi--Fredholm for $\lambda \in {\Bbb C}\setminus 
{\Bbb R}.$ Thus $\dim E_\gl(A)$ is locally 
constant on $\C\setminus \R$ (see \cite{Kat:PTL}) and therefore 
$\dim E_{\pm\gl}(A)=\dim E_{\pm i}(A)$ 
for $\gl\in \C_+$. For another proof see Corollary \plref{cor1.1}
below.

There are situations in which it is clear that the formal
defect spaces $\ce_\gl(S)$ and the defect spaces $E_\gl(S)$ are isomorphic.
This is, for instance, the case if $\ch(x)$ is invertible for almost all $x\in
I$. In general, the analogue of \eqref{SDM-G1.14} for the dimensions of the formal 
defect subspaces $\ce_{\gl}(S)$ holds.
However, this is less trivial. The only proof we know of so far is due to Kogan and
Rofe--Beketov \cite[Sec. 2]{KogRof:SIS}. It
uses methods from complex analysis and is rather technical.
Here we can give a very simple proof of this fact which is based on
the regularity Theorem \plref{regthm}. Namely, the regularity Theorem allows 
to show a simple relation between the deficiency indices and the formal 
deficiency indices:
     \begin{prop}\label{ML-S2.19} 
Let $S$ be a general symmetric system \eqref{G1.2} 
on an interval $I\subset\R$. Then for $\gl\in\C$ we have
\begin{equation}
   \dim\ce_\gl(S)=\dim E_\gl(S)+n-\rank S.
\end{equation}
In particular, if the system is definite then $\dim\ce_\gl(S)=\dim E_\gl(S)$. 
\end{prop}
\begin{proof}
Consider $\tilde f\in E_\gl(S)$. This means $\{\tilde f,\gl\tilde f\}\in \Smax$
and in view of Theorem \plref{regthm} there exists
$f\in\tilde f,f\in \AC(I,\C^n)\cap \cLH$ 
such that $Jf'+Bf=\gl\ch f$.
Thus $f\in \ce_{\lambda }(S)$. 
This shows that the quotient map $\pi:\ce_\gl(S)\to E_\gl(S)$ is surjective.

Next let $\{f,\gl f\}\in \ker \pi$. This means that 
$Jf'+Bf=\gl\ch f$ and $\tilde f=0$. Thus $\ch f=0$. Hence $\ker \pi$
consists of the solutions of $Jf'+Bf=0, \ch f=0$. This space is isomorphic
to $\ker\Phi$ (cf. Subsections \plref{sec21}, \plref{sec22}) and
hence $\dim\ker\pi=\dim\ker\Phi=n-\rank S$ and we reach the conclusion.
\end{proof}

The following result was proved by Kogan and Rofe--Beketov for the 
half--line \cite[Theorem 2.1]{KogRof:SIS} and for systems
on the line which are definite on both half--lines
$\R_\pm$ \cite[Corollary 2.2]{KogRof:SIS}. For general non--definite
systems it seems to be new. 
     \begin{prop}\label{S1.11}
Let $S$ be a general symmetric system \eqref{G1.2} on
an interval $I\subset \R$. Then 
\[\dim \ce_{\pm \lambda }(S)=\dim \ce_{\pm i}(S)=:\cn_{\pm 
}(S),\quad\text{for}\quad \gl\in\C_+.
\]
\end{prop}
\begin{proof} This follows immediately from
\eqref{SDM-G1.14} and Proposition \plref{ML-S2.19}.
    \end{proof}

For completeness we note the case of a finite 
interval:

     \begin{prop}\label{S1.9} Let $I=[a,b]$ be a finite interval and $S$
the symmetric linear relation obtained from the first order
system \eqref{G1.2}, where $B,\ch\in \cl^1(a,b)$. 
Then $\cn_\pm(\cs)=n$ and $N_\pm(S)=\rank (S)$. 

In particular $\cn_{\pm}({\cs})=N_{\pm}(S)=n$ if the system $S$ is definite.
      \end{prop}
   \begin{proof} It is clear that the differential 
equation \eqref{G1.2} 
has $n$ linear independent solutions. Hence 
$\cn_\pm(\cs)=n$. From
Proposition \plref{S1.6} we infer that 
$\widetilde{\fsys\xi}$ is nonzero
if and only if $\xi\in\im \Mmat$. This implies 
$N_\pm(S)=\rank \Mmat$.
   \end{proof}
\subsection{Definite systems and von Neumann formula}\label{sec25}

We start with the following generalization of the von Neumann formula.

  \begin{prop}\label{Mark-S1.10} 
Let $A$ be a closed symmetric linear relation in the
Hilbert space ${\frak H}.$ Then for each
pair $\{\gl _1,\gl _2\}\in {\mathbb C}_+\times 
{\mathbb C}_-
$ we have the following 
direct sum decomposition
    \begin{equation}\label{G1.34}
A^*=A\dotplus {\hat E}_{\gl _1}\dotplus {\hat E}_{\gl_2},
\qquad {\hat E}_{\gl}=
\bigsetdef{\{f,\gl f\}}{f\in E_{\gl }}.
   \end{equation}
  \end{prop}
$\dotplus$ denotes a (non--orthogonal) direct sum of 
vector spaces.

\begin{proof}
We put $A_{\gl_1}:=A\dotplus{\hat E}_{\gl_1}.$ 
It is clear that
$A_{\gl_1}$ is a closed dissipative extension of $A,$ 
that is 
$A\subset A_{\gl_1}\subset A^*$ and $\Im(f,g)\ge 0$ 
for any
$\{f,g\}\in A_{\gl _1}.$ In fact, we show that 
$A_{\gl _1}$ is a maximal 
dissipative relation in ${\frak H}.$ 
To prove this fact it suffices to check that
${\ovl \gl}_1\in\rho(A_{\gl_1}),$ where $\rho(A_{\gl_1})$
denotes the resolvent set of $A_{\gl_1}$. 

For a dissipative linear relation $T$ and 
$\mu =\alpha -i\beta \in {\mathbb C}_-$ one has for 
$\{f,g\}\in T$
   \begin{equation*}
\|g-\mu f\|^2=\|g-\alpha f\|^2+2\beta 
\Im(f,g)+\beta^2\|f\|^2\ge
\beta ^2\|f\|^2.
  \end{equation*}
Hence $A_{\gl_1}-\ovl{\gl_1}I$ is injective with 
closed range and thus it suffices
to verify that $\im(A_{\gl _1}-{\ovl {\gl}}_1I)$ is 
dense in ${\frak H}$.

Let $\varphi $ be orthogonal to $\im(A_{\gl _1}-
{\ovl{\gl}}_1I)$,
that is 
   \begin{equation}
\scalar{g-{\ovl 
{\gl}}_1f}{\varphi}=0\quad\text{for}\quad \{f,g\}\in 
A_{\gl _1}.
\label{ML-G2.47}
    \end{equation}
In particular, we have for $\{f,g\}\in A$
\[ 
\scalar{g}{\varphi}=\scalar{\ovl{\gl_1}f}{\varphi}=
\scalar{f}{\gl_1\varphi}.
\]
Hence $\varphi \in E_{\gl _1}$ and 
$\{\varphi ,\gl _1\varphi \}\in {\hat E}_{\gl_1}.$
>From the latter and \eqref{ML-G2.47}
we infer $0=\scalar{\lambda _1\varphi -{\overline 
{\lambda }}_1\varphi}{\varphi}=
-2i\Im\gl _1\|\varphi\|^2.$ Hence  $\varphi =0$. 
Summing up, we have proved that ${\mathbb C}_-\subset 
\rho (A_{\gl_1})$ and hence $A_{\gl_1}$
is maximal dissipative.

On the other hand for each proper extension 
${\widetilde A},  A\subset {\widetilde A}\subset A^*$ 
the inclusion 
$\mu \in \rho ({\widetilde A})$ is equivalent to the 
fact that 
${\widetilde A}$ is transversal 
\footnote{Two proper extensions $A_1$ and 
$A_2$ of $A$ are called transversal if $A_1\cap A_2=A$ 
and $A_1+A_2=A^*.$}
to $A_{\mu }:=A\dotplus{\hat E}_{\mu}$ (see 
\cite{Mal:FGR}). Hence $A_{\lambda _1}$ and 
$A_{\gl _2}$ are transversal and this is equivalent to 
the direct sum
decomposition (\ref{G1.34}).
   \end{proof}
Now we can give the second proof of \eqref{SDM-G1.14}.
    \begin{cor}[{\cite{AkhGla:TLO}, 
\cite{Naj:LDO}}]\label{cor1.1}
With the previous notations we have for all 
$\gl\in\C_+$
   \begin{equation}\label{G1.35}
\dim E_{\pm \lambda }=\dim E_{\pm i}.
    \end{equation}
   \end{cor}
    \begin{proof}
Let $\gl _2=-i.$ It follows from (\ref{G1.34}), that 
for each 
$\gl_1\in \C_+$ 
   \begin{equation*} 
\dim E_{\lambda _1}=\dim A^*/(A\dotplus \hat E_{-i}).
    \end{equation*}
     \end{proof}
\begin{remark}\label{rem2.5}
1. Formula (\ref{G1.34}) with $\gl_2={\ovl{\gl}}_1$ is 
well--known
\cite{Orc:CDE}, \cite{Cod:ETF}, \cite{Ben:SRH}. For 
$\gl_1=i={\ovl{\gl}}_2$ the direct sum
(\ref{G1.34}) is orthogonal 
        \begin{equation}\label{G1.36}
A^*=A\oplus \hat E_i\oplus \hat E_{-i}.
        \end{equation}

2. The maximal dissipativity of the linear relation 
$A_{\gl}$ with
$\gl \in {\mathbb C}_+$ is well--known. We presented 
the 
proof for the sake of 
completeness. Note, however, that our proof of this 
fact as well as the proof 
of the well--known Corollary \ref{cor1.1} is simpler 
and shorter than the known ones.
   \end{remark}

We continue in noting a simple lemma which is a 
generalization of a well--known result
(cf. \cite{Naj:LDO}) on symmetric operators to the 
case of symmetric linear relations.

\begin{lemma}\label{le3.1}
Let $A$ be a closed symmetric linear relation in the 
Hilbert space $\frak H$ and
$\ker(A-aI)=\{0\}$ for some $a\in {\mathbb R}$. Then
   \begin{equation}\label{3.3}
\dim E_a(A)=\dim \ker (A^*-aI)\le N_{\pm }(A).
   \end{equation}
  \end{lemma}
         \begin{proof}
Similar to the proof of Proposition \plref{Mark-S1.10} 
we put
   \begin{equation*}
{\widetilde A}_a:=A\dotplus \hat E_a(A) ,\qquad
\hat E_a(A):=\bigsetdef{\{f,af\}}{f\in E_a(A)}.
   \end{equation*}
It is clear that ${\widetilde A}_a$ is a symmetric 
extension of $A$ and the subspaces $A$ and $\hat E_a(A)$ are 
linearly independent 
since $\ker(A-aI)=\{0\}.$ Therefore
$\dim({\widetilde A}_a/A)=\dim E_a(A)$. On the
other hand the von Neumann formula for linear 
relations \eqref{G1.36}
yields $\dim({\widetilde A}_a/A)\le \min(N_+,N_-).$ 
Combining these relations we obtain \eqref{3.3}.
   \end{proof}
We return to the discussion of the relation $S=S(J,B,\ch)$.
Denote by $\Exts(\cs)$ and $\Exts(S)$ the set of 
closed symmetric extensions 
of $\csmin$ and $\Smin$ respectively:
    \begin{equation}\begin{split}
          \Exts(\cs)&=\bigsetdef{\widetilde 
{\cs}}{\csmin\subset {\widetilde {\cs}}
     \subset \csmax, \tilde\cs\text{ is closed and 
symmetric}},\\
          \Exts(S)&=\bigsetdef{\widetilde 
S}{\Smin\subset 
{\widetilde S}\subset \Smax, \tilde S\text{ is closed 
and symmetric}}.
		    \end{split}
     \end{equation}
   \begin{prop}\label{p3.1} 
Assume that the system \eqref{G1.2} is definite on 
${\mathbb R}_+$. Then:
   \begin{thmenum}
\item The quotient map $\pi$ maps $E_{\lambda}(S)$ 
isomorphically onto
$\ce_{\gl}(S)$ for each ${\gl}\in 
{\mathbb 
C}$ and consequently
$\cn_{\pm }(S)=N_{\pm }(S).$

\item  For each $\xi \in {\Bbb C}^n$ and each $a\in 
[0,\infty]$ 
there exists $\{f,g\}\in {\cs}^*$ with compact support 
such that 
$f(a)=\xi.$ If $a>0$ and the system $S$ is definite 
both on $[0, a]$ and $[a, \infty]$, then $\{f,g\}$ can
be chosen such that $\{f,g\}\in {\csmin}$.

\item  If $\{{\widetilde f},{\widetilde g}\}\in S$ 
with representatives
$\{f,g\}\in {\csmax}$ then $f(0)=0$, that is 
      \begin{equation*}
\cd({\cs})\subset \bigsetdef{f\in \AC({\R}_+,{\C}^n)}{f(0)=0}.
    \end{equation*}

\item  The quotient map $\pi _2:=\pi \oplus \pi \ 
\text {maps}\ {\cs}^*$ and 
$\cs$ isomorphically onto $\Smax$ and $S$ 
respectively.

\item  For each pair $\{{\gl}_1,{\gl}_2\}\in {\mathbb 
C}_+\times {\mathbb C}_-$ 
the following analogue of the von Neumann formula holds  
true
        \begin{equation}\label{1.41}
\csmax=\cs\dotplus \hat \ce_{\gl_1}(\cs)\dotplus
\hat \ce_{\gl_2}(\cs), \quad
   \hat \ce_\gl(\cs):=\bigsetdef{\{f,\gl f\}}{ f\in 
\ce_{\gl }}.
        \end{equation}
For ${\gl}_1=i={\ovl{\gl}}_2$ the direct sum 
decomposition \eqref{1.41} is
orthogonal.

\item  The quotient map $\pi _2:=\pi \oplus \pi $ 
induces a bijective 
correspondence between the sets $\Exts(\cs)$ and 
$\Exts(S)$. Moreover, 
${\tilde {\cs}}$ is self--adjoint iff ${\tilde S}$ 
is self--adjoint.

\item For each $a\in {\mathbb R}$ the following 
inequality holds
   \begin{equation}\label{3.5}
\cn_{\pm }({\cs})\ge \dim \ce_a({\cs})=N_a(S).
   \end{equation}
\end{thmenum}
   \end{prop}
    \begin{proof}
(1) has been established in the proof of Proposition 
\plref{ML-S2.19}.  

(2) is implied by Proposition \plref{S1.8} (3) since 
$\im\Phi={\Bbb C^n}$ .  

(3) is a special case of Proposition \plref{S1.8} (4) 
since $\ker \Phi={0}$ . 

(4) Injectivity of the map $\pi _2:{\cs}^*\to \Smax$ 
follows again from the 
assumption that $S$ is definite. Indeed, let 
$\{{\widetilde f},{\widetilde g}\}\in \Smax,  
f_k\in {\widetilde f}, \{f_k, g\}\in {\cs}^*, k=1,2.$ 
Then $f:=f_1-f_2$ satisfies the homogeneous equation
$Jf'+ Bf=0$, that is $f\in E_0({\cs})$. Since 
$f_1,f_2\in\tilde f$
we have $\ch f=0$ and therefore 
$\int _If^*(x){\ch}(x)f(x)dx=0.$
Since $S$ is definite the latter implies $f=0$.

Surjectivity has been established in Proposition 
\plref{S1.8}.

(5) is a consequence of (1),(4) and Proposition 
\plref{Mark-S1.10}.

(6) W.l.o.g. we may assume $N_+\le N_-$. By 
definiteness we then have
$\cn_+=N_+\le N_-=\cn_-.$ It follows from the von Neumann
formula \eqref{1.41} with ${\gl}_1={\ovl {\gl}}_2=i$ 
that each symmetric 
extension ${\widetilde {\cs}}\supset \cs$ is 
given by the second Neumann formula
    \begin{equation}\label{1.42}
{\tilde {\cs}}=\cs\dotplus \bigl\{(I+V)\ce'_i,\ 
i(I-V)\ce'_i\bigr\}, 
    \end{equation}
where $\ce'_i\subset \ce_i$ is a linear subspace and $V$ 
is an isometric operator from $\ce'_i$ onto
$V\ce'_i\subset  \ce_{-i}.$ The corresponding symmetric 
extension 
${\tilde S}\supset S$ is given by 
    \begin{equation}\label{1.43}
{\tilde S}=S\dotplus \bigl\{(I+V)E'_i,\  
i(I-V)E'_i\}, \qquad  E'_i:=\pi 
\ce'_i\subset E_i.
     \end{equation}  
It is clear from (1) that this establishes the 
asserted bijective correspondence. Cf. also
Proposition \plref{ML-S2.15}

(7) We know from Proposition \plref{S1.8} (1) that 
$\ker(S-aI)=\{0\}.$ 
Lemma \plref{le3.1} and (1) now imply
\begin{equation}
  \cn_\pm({\cs})=N_\pm(S)\ge \dim\ker(\Smax-aI)=\dim 
\ce_a({\cs}).
\end{equation}
      \end{proof}     

\section{Essential self--adjointness on the line. First 
approach.}\label{sec3}

\subsection{Preliminaries and a first criterion for essential self--adjointness}\label{sec31}

In this section we study the system \eqref{G1.2} on 
the real line and discuss essential self--adjointness. For the 
moment let $I\subset \R$
be an interval and $S=S(J,B,\ch)$ be the symmetric linear 
relation of the first order
system \eqref{G1.2}. 
Let 
\begin{equation}
   \gl_j(x):=\max_{V\subset \C^n, \dim
   V=j-1}\min\bigsetdef{\scalar{\ch(x)\xi}{\xi}}{%
          \xi\perp V, \|\xi\|=1}
\label{ML-G3.1}
\end{equation}
be the $j-$th eigenvalue of $\ch(x)$. 
Furthermore, we put
\begin{equation}
    c(x):=\begin{cases}  \|\ch(x)^{-1/2}J(x)\ch(x)^{-1/2}\|,&
    \det(\ch(x))\not=0,\\
          \infty,&\textup{otherwise}.
          \end{cases}
\label{G2.3}
\end{equation}
We have estimates
\begin{equation}
   \frac{1}{c(x)}\le \|\ch(x)^{1/2}J(x)^{-1}\ch(x)^{1/2}\|
               \le \|J(x)^{-1}\| \gl_n(x),
    \label{ML-G3.3a}
\end{equation}
and, if $\det{\ch(x)}\not=0$,
\begin{equation}
     c(x)\le \|J(x)\|\|\ch(x)^{-1}\|=\frac{\|J(x)\|}{\gl_1(x)}.
     \label{ML-G3.3b}
\end{equation}
Thus we have for all $x\in \R$
\begin{equation}
  \frac{\gl_1(x)}{\|J(x)\|}\le \frac{1}{c(x)}\le \|J(x)^{-1}\| \gl_n(x).
  \label{ML-G3.3}
\end{equation}
In view of \eqref{ML-G3.3a} the function $\frac{1}{c(x)}$ 
is locally integrable. The significance of $c(x)$ stems from the fact that
if $\det(\ch(x))\not=0$ then
for $\xi\in\C^n$ we have the estimate
\begin{equation}\begin{split}
      \xi^* J(x)^* \ch(x)^{-1} J(x) \xi&=\|\ch(x)^{-1/2} 
J(x)\xi\|^2\le c(x)^2 \|\ch(x)^{1/2}\xi\|^2\\
       &= c(x)^2 \xi^* \ch(x) \xi.
		\end{split}
\label{ML-G3.5}
\end{equation}

\begin{lemma}\label{ML-S2}
Let $f\in L^1_{\loc}(\R), f(x)\ge 0$, be a non--negative
locally integrable function.
Assume in addition that
\[\int_0^\infty f(x) dx=+\infty.\]
Then for $n\in\N$ there exists an absolute continuous 
function
$\chi_n\in\AC(\R)$ with the properties
\begin{thmenum}
\item $\chi_n(x)=1, x\le n,$
\item $\chi_n(x)=0, x\ge x_n$, for some $x_n$,
\item $\chi_n'\in L^\infty(\R),$
\item $|\chi_n'(x)| \le \frac 1n f(x)$, for all 
$x\in\R$.
\end{thmenum}
\end{lemma}
\begin{proof}
Fix $n\in\N$. By B. Levy's theorem on monotone 
convergence we have
\[ \lim_{C\to+\infty}\int_0^\infty \frac 1n 
\min(C,f(x))dx=+\infty,\]
and thus we may choose $C>0$ such that 
\[   \int_n^\infty  \min(C,\frac 1n f(x))dx \ge 2.\]
Now choose $N$ large enough such that
\[ K_n:=\int_n^N  \min(C,\frac 1n f(x)) dx \ge 1\]
and put
\[\chi_n(x):=1- \frac {1}{K_n} 
\int_{\min(n,x)}^{\min(N,x)}\min(C,\frac 1n f(s))ds.\]
$\chi_n$ has the desired properties with $x_n=N$.
   \end{proof}
      \begin{theorem}\label{S2.2}
Let $\Smin=\Smin(J,B,\ch)$ be a first order
system \eqref{G1.2} on the interval $I$.
\begin{thmenum}
\item Let $I=\R$ and assume that
\begin{equation}
 \pm\int_0^{\pm\infty} \frac{1}{c(x)} dx=+\infty.\label{ML-G1.6}
\end{equation}
Then $\Smin$ is essentially self--adjoint, i.e. 
$S=\Smax$.
\item Let $I=\R_+$ and assume that
\begin{equation}
\int_0^{\infty} \frac{1}{c(x)} dx=+\infty.\label{ML-G1.6a}
\end{equation}
Then for $\{\tilde f,\tilde g\}\in \Smax$ there exists 
a sequence
$\{f_n,g_n\}\in \csmax$ such that $\tilde f_n\to 
\tilde f, \tilde g_n\to \tilde
g$ in $\LHv{\R_+}$ and $\supp f_n, \supp g_n\subset 
[0,\infty)$ compact.
Moreover, for $\{\tilde f_j,\tilde g_j\}\in \Smax, 
j=1,2$ one has
\begin{equation}
             \scalar{\tilde f_1}{\tilde g_2}-
\scalar{\tilde g_1}{\tilde
             f_2}=-f_1(0)^*Jf_2(0).
   \label{G2.8}
\end{equation}
\end{thmenum}
         \end{theorem}
      \begin{remark}\label{ML-S3.3} Note that the condition
\eqref{ML-G1.6} 
(resp. \eqref{ML-G1.6a})
implies that for each $R>0$ there exist subsets $K_\pm\subset \R_\pm\setminus
[-R,R]$ (resp. $K_+\subset \R_+\setminus[0,R]$) of positive Lebesgue measure 
such that the Hamiltonian $\ch$ is positive definite 
on $K_{\pm}$ (resp. $K_+$).
In particular, the corresponding system $S$ is 
definite on $\mathbb R$ (resp. on $\mathbb R_+$).  
         \end{remark}
\begin{proof} (1) According to Lemma \plref{ML-S2} let 
$\chi_n$ be absolutely continuous
with bounded derivative,
\[ \chi_n(x)=\begin{cases} 1,& |x|\le n,\\
                            0, & |x|\ge x_n,
              \end{cases}
\]
and
\[   |\chi_n'(x)| \le \frac{1}{n c(x)}.\]
For $\{\tilde f,\tilde g\}\in \Smax$ we choose, 
according to Proposition
\plref{S1.8}, representatives $\{f,g\}\in\csmax$ and 
put
\[
     f_n:= \chi_n f.
\]
Since $\chi_n'$ vanishes if $\ch(x)$ is not invertible 
the function
$\chi_n' \ch(x)^{-1} J(x)f$ is well--defined. Moreover
\begin{align*}
    \|\chi_n' \ch^{-1} J f\|_{\LHv{\R}}^2&\le \int_\R 
|\chi_n'(x)|^2
                f(x)^*J(x)^* \ch(x)^{-1}J(x) f(x) dx\\
     &\le \sup_{x\in \R}( \chi_n'(x) c(x))^2 
\|f\|_{\LHv{\R}}^2\\
     &\le \frac{1}{n^2} \|f\|_{\LHv{\R}}^2,
\end{align*}
hence $\chi_n' \ch(x)^{-1} J(x)f$ lies in $\cLHv{\R}$ and 
it converges to $0$ in $\cLHv{\R}$.
Finally, we calculate
\begin{align*}
     Jf_n' +Bf_n &= \chi_n (J f' + Bf)+ \chi_n' Jf\\
                &= \ch(\chi_n g+\chi_n' \ch^{-1} Jf)\\
                &=:\ch g_n.
\end{align*}
Thus $\{f_n,g_n\}\in \cs$ and 
$\lim\limits_{n\to\infty}\{\tilde f_n,\tilde
g_n\}=\{\tilde f,\tilde g\}$ and the claim is proved.

The proof of (2) proceeds along the same lines with 
minor modifications.

\eqref{G2.8} follows from integration by parts if 
$f_2,g_2$ have compact support.
To prove it in general we consider $f_{2,n}=\chi_nf_2$ 
and $g_{2,n}=\chi_n g_2+
\chi'\ch^{-1} J f_2$. Then \eqref{G2.8} holds true for 
$\{\tilde f_1,\tilde
g_1\}$
and $\{\tilde f_{2,n},\tilde g_{2,n}\}$. Noting that 
$f_{2,n}(0)=f_2(0)$ is independent
of $n$ we obtain the result by taking the limit as 
$n\to \infty$.
    \end{proof}
    \begin{remark}\label{rem2.2} 
 \eqref{ML-G1.6} is not necessary for $\Smin$ to be essentially
self--adjoint. Namely, in the situation of Example 
\plref{Mark-S1.11} 2. we have
$\frac{1}{c(x)}=0$. But there certainly exist $V$ (e.g. $V\in 
L^\infty(\R)$) such 
that the Schr\"odinger operator $-\frac{d^2}{dx^2}+V$ 
and hence the Hamiltonian in Example \plref{Mark-S1.11} 2. 
are essentially self--adjoint. 

See also Example \plref{ML-S5.33} for a counterexample with a 
nonsingular Hamiltonian ${\ch}$.
\end{remark}
\begin{cor}\label{S2.3} If $J=J(0)$ is constant then the condition
\eqref{ML-G1.6} (resp. \eqref{ML-G1.6a}) is implied by
\begin{equation}
    \pm\int_0^{\pm\infty} \gl_1(x)dx=\infty \qquad 
(\textrm{resp.}\;
    \int_0^\infty \gl_1(x)dx=\infty).\label{G2.9}
\end{equation}
Hence (for $J=J(0)$) \eqref{G2.9} implies the conclusions in Theorem 
\plref{S2.2}.
\end{cor}
\begin{proof} This follows immediately from the estimate
\eqref{ML-G3.3}.
\end{proof}
     \begin{remark}\label{rem2.2'}
It is clear that $\Smin$ is essentially self-adjoint iff for each
$f \in \cd(\csmax)$ the following limit exists:
   \begin{equation}
\lim\limits_{x\to \pm \infty }f(x)^*J(x)f(x)=0.
     \end{equation}

Condition \eqref{G2.9} yields a stronger conclusion about $\cd(\csmax).$
In order to explain it we denote by $\AC_0({\mathbb R},{\mathbb C}^n)$ the 
set of those $f\in \AC({\mathbb R},{\mathbb C}^n)$
such that there exist sequences $x^{\pm }_n\to \pm \infty$ 
with 
$\lim\limits_{n\to \infty}\scalar{f(x^{\pm}_n)}{f(x^{\pm}_n)}=0.$ 

It is clear that under condition \eqref{G2.9} 
$\cd(\csmax)\subset \AC_0({\mathbb R},{\mathbb C}^n).$ 
The converse assertion is also true if $B=0$ and
$\ch(x)=\diag\bigl(\gl_1(x),...,\gl_n(x)\bigr)$ is a diagonal matrix
with eigenvalues $\gl_1(x)\le ...\le \gl_n(x).$ Indeed, if
$\gl _1\in L^1({\Bbb R})$ then 
   \begin{equation*}
f:=\col(1,0,...,0)\in \cd(\csmax)\ \text{but}\ \ 
\scalar{f(x)}{f(x)}=1\ \text{and}\ \ 
f{\notin}\AC_0({\mathbb R},{\mathbb C}^n).
    \end{equation*}
Probably \eqref{G2.9} is equivalent to the inclusion
$\cd(\csmax)\subset \AC_0({\mathbb R},{\mathbb C}^n)$ for an arbitrary 
$\Smax(J,B,\ch)$ with constant $J=J(0).$

   \end{remark}

However \eqref{G2.9} is weaker than \eqref{ML-G1.6} 
as the following example shows:

\begin{example}\label{ex3.6} 
Let
  \begin{equation*}
\ch(x)=\diag(\gl_1(x),\gl_2(x)), \quad  
J=\begin{pmatrix} 0&1\\
-1&0
\end{pmatrix},
\end{equation*}
where
  \begin{equation*}
\gl_1(x)=(|x|+2)^{-1}\ln^{-2}(|x|+2) , \quad 
\gl_2(x)=(2+|x|)^{-1}.
   \end{equation*}
Then $\gl_1(x)\in L^1(\R)$, but
 \begin{equation*}
(\gl_1\gl_2)^{1/2}\not\in L^1(\R_\pm).
  \end{equation*}
Since $\pm i (\gl_1\gl_2)^{1/2}$ are the 
eigenvalues of
$\ch^{1/2}J^{-1}\ch^{1/2}$ we infer that 
$\frac{1}{c(x)}=(\gl_1(x)\gl_2(x))^{1/2}$.
Hence \eqref{ML-G1.6} is satisfied but $\gl_1\in L^1(\R)$.

Besides, setting 
     \begin{equation*}
f=\col\bigl(\ln^{1/4}(2+|x|), 0\bigr), \quad    
g=\col\bigl(0, -4^{-1}\sgn(x)\ln^{-3/4}(2+|x|)\bigr)\ \in \cLHv{\R} 
     \end{equation*}
one gets
$\{f,g\}\in \Smax$ but $\scalar{f(x)}{f(x)}=\ln^{1/2}(2+|x|)\to \infty \ 
\text{as}\ x\to \pm \infty\ \text{and}\ f{\notin}
\AC_0({\mathbb R},{\mathbb C}^2).$
     \end{example}
    \subsection{The case of a symmetric operator}\label{sec32}
For completeness we briefly comment on the case
that the system \eqref{G1.2} defines a symmetric 
linear operator
containing at least the $C^1$--functions 
with compact support
in its domain. Namely, let $J,B,\ch$ be as in \eqref{G1.0}
and assume in addition that $\ch(x)$ is invertible for all
$x\in I$ and that $\ch(x)^{-1}, B^*\ch^{-1}B$ is locally integrable.
In this case each class $\tilde f\in\LH$ contains at most
one continuous representative. In particular $\AC_{\comp}(I,\C^n)$
may be viewed as a subset of $\LH$.

Then we consider the differential operator
  \begin{equation}
       L:=\ch^{-1}\bigl(J\frac{d}{dx} +B\bigr).\label{ML-G1.1}
  \end{equation}
$L$ maps $\cd(L):=C^1_{\comp}(I,\C^n)$ into $\LH$. Namely, if
$f\in C^1_{\comp}(I,\C^n), K:=\supp(f)\subset I$ then we estimate
   \begin{equation}\begin{split}
      \|\ch^{-1}Jf'\|_{\ch}^2&=\bigl|\int_Kf'(x)^*J(x)\ch(x)^{-1}J(x)f'(x)dx
      \bigr|\\
      &\quad\le \sup_{x\in K}\|J(x)f'(x)\|^2 \int_K \|\ch(x)^{-1}\|dx<\infty\\
      \|\ch^{-1}Bf\|_{\ch}^2&=\int_Kf(x)^*B(x)^*\ch(x)^{-1}B(x)f(x)dx\\
      & \quad\le
      \sup_{x\in K} \|f(x)\|^2 \int_K \|B(x)^*\ch(x)^{-1} B(x)\|dx<\infty.
		\end{split}
   \end{equation}
$L$ is formally symmetric and in view of the regularity Theorem \plref{regthm}
the domain $\cd(L_{\max})$ of $L_{\max}:=L^*$ lies in $\AC(I,\C^n)$. 
Furthermore, for $f,g\in\cd(L_{\max})$ and $\ga<\gb$ we have
\begin{equation}\begin{split}
    &\int_\ga^\beta (L_{\max}f)(x)^* \ch(x) g(x) dx-
      \int_\ga^\beta f(x)^* \ch(x)(L_{\max}g)(x) dx\\
     &=  -f(\gb)^* J(\gb) g(\gb)+
         f(\ga)^* J(\ga) g(\ga).
		\end{split}
\label{ML-G3.26}
      \end{equation}
In contrast to general first order systems the domain of $L_{\max}$ is
localizable in the following sense: $\cinfz{I,\C^n}$ is dense in 
\[\cd_{\comp}(L_{\max}):=
\bigsetdef{f\in\cd(L_{\max})}{\supp(f) \;\text{is 
compact}}\]
with respect to the graph norm of $L$. Namely, from \eqref{ML-G3.26} we infer 
that for 
$f\in\cd_{\comp}(L_{\max}), g\in\cd(L_{\max})$
we have 
$\scalar{L_{\max}f}{g}=\scalar{f}{L_{\max}g}$, i.e.
$f\in \cd(L_{\max}^*)=\cd(L^{**})=\cd(\ovl{L})$.

Summing up one arrives at the following result.
      \begin{theorem}\label{ML-S1.3}
Let $I=\R$ and let $\ch(x)$ be invertible for $x\in \R.$ Assume also 
that $\ch^{-1}$
and $B^*\ch^{-1}B$\ are locally integrable and  \eqref{ML-G1.6} holds. 
Then the operator $L$ is essentially self--adjoint on
$\cinfz{\R,\C^n}\subset L^2_\ch(\R,\C^n)$. 
   \end{theorem}

    \begin{cor}\label{cor3.14}
In the framework of Theorem \plref{ML-S1.3} assume that $J(x)$ is bounded
on $\R$ and that there exists a $\delta>0$ such that 
$\ch(x)\ge\delta>0$ for $x\in\R$. Moreover, assume 
$B\in L^2_{\loc}(\R,\Mat(n,\C))$. 
Then $L$ is essentially self--adjoint on
$\AC_{\comp}(\R,\C^n)\subset L^2_\ch(\R,\C^{2n})$.
     \end{cor}
\begin{proof} $\ch(x)\ge\delta>0 $ implies that 
$\ch(x)^{-1/2}\le \delta^{-1/2}.$
Hence \eqref{ML-G1.6} holds since 
$\|\ch(x)^{-1/2}J(x)\ch(x)^{-1/2}\| \le C\delta^{-1}.$ 
Moreover, from $\ch(x)\ge \delta>0$ we infer that $\ch^{-1}$
is bounded and hence $B^*\ch(x)^{-1}B\le \frac 1\delta B^*B\in L^1_{\loc}$.
Hence Theorem \plref{ML-S1.3} applies.\end{proof}

      \begin{cor}\label{cor3.15}
In the framework of Theorem \plref{ML-S1.3}
let $J(x)$ be bounded on $\R$ and let
$\ch(x)=T^*(x)\ch_1(x)T(x)$ such that
\begin{thmenum}
\item $T(x)$ and $\ch_1(x)$ are continuous on $\R$,
\item $\ch_1(x)\ge \delta >0, \quad x\in \R$,
\item $T^*(x)Q(x)T(x)=J(x)$, where $Q$ is continuous and bounded.
\item $B\in L^2_{\loc}(\R,\Mat(n,\C)$.
\end{thmenum}
Then $L$ is essentially self--adjoint on
$\AC_{\comp}(\R,\C^n)\subset L^2_\ch(\R,\C^{2n})$.
      \end{cor}
\begin{proof} Since $T^*QT=J$ it is clear that $T(x)$ 
and $Q(x)$ are invertible
for all $x$. Furthermore, 
$\ch(x)=T^*(x)\ch_1(x)T(x)\ge \delta T^*(x)T(x),$ 
hence
$\ch(x)>0$ for all $x$.  Setting $K(x):=\ch^{-
1/2}(x)T^*(x)$, one has
$\|K(x)\|\le \delta ^{-1/2}$
and thus
\begin{align*}
c(x) &\le \|\ch^{-1/2}(x)J(x)\ch^{-1/2}(x)\|\\
     &\le \|K(x)Q(x)K^*(x)\|\\
     &\le \delta^{-1}\|Q(x)\|\le c\delta ^{-1}
\end{align*}
since $Q(x)$ is bounded on $\R$.
Hence \eqref{ML-G1.6} is fulfilled and we reach the 
conclusion.
\end{proof}
     \begin{remark}\label{rem3.16} 
1.  If $\ch(x)$ is invertible for almost all $x\in \R$ then by Theorem 
\plref{S2.2} the operator $L_{\min}$ defined by \eqref{ML-G1.1} on 
     \begin{equation*}
\cd(L_{\min})=\cd(\Smin)=\{f\in \AC_{\comp}(\R,\C^n)| \ \ Lf\in L^2_{\ch}(\R)\}
   \end{equation*}
is essentially self--adjoint under the only condition \eqref{ML-G1.6}. 

However, we cannot conclude the essential self--adjointness of $L$
on $\cinfz{\R,\C^n}$ without additional assumptions 
(like in Theorem \plref{ML-S1.3}) since
in general $C^1_{\comp}(\R,\C^n)$ is not contained in $\cd(L_{\min})$.

2. Corollary \plref{cor3.14} and  Corollary \plref{cor3.15} have been 
obtained by L. Sakhnovich {\cite{Sak:DIS}} under the additional assumptions
$B=0$ and $J=J(0)$ constant.

3. In {\cite[Proposition 2.1]{LesMal:ISP}} we established self--adjointness 
of the operator $L$ with ${\ch}=I$ and $J(x)=J(0)$ 
being constant. This fact is well--known. 
It is contained, e.g., as a very special case 
in a result due to Levitan and Otelbaev {\cite[Theorem 2]{LevOte:SAC}}.

Note however that the proof of Proposition 2.1 from {\cite{LesMal:ISP}} 
remains valid  if $J(x)$ is nonconstant and bounded on $\R.$ 
Corollary \plref{cor3.14} 
is reduced to this result via the gauge transformation \eqref{2.8} 
with $U={\ch}^{-1/2}.$
\end{remark}

\section{Essential self--adjointness on the line. Second 
approach.}\label{sec4}

In this section we present a second proof of the essential
self--adjointness of the operator $L$ from Subsection \plref{sec32}.
This second proof uses the hyperbolic equation method
(cf. \cite{Ber:EES},\cite{Che:ESP}).

If the coefficients of $L$ are smooth then this
method even proves the 
essential self--adjointness
of all powers $L^n(n\in \Z_+)$ of the operator $L$
\cite{Che:ESP}.

We recall some definitions and results.
Let $H$ be a densely defined operator in a Hilbert 
space $\frak H$.
Recall 
that a vector function $u:[0,\infty )\to \frak H$ is 
called a strong
solution 
of the equation
    \begin{equation} \label{2.1}
\frac{du}{dt}(t)+Hu(t)=0, \qquad  t\in (0,\infty ),
      \end{equation}
if $u$ is strongly differentiable, $u(t)\in \mathcal 
D(H)$ for each 
$t\in (0,\infty )$ and \eqref{2.1} is satisfied for 
each $t\in (0,\infty )$.

Our second proof of the essential self--adjointness is 
based on the following result due to 
Berezanskii-Povzner (cf. also \cite{Che:ESP}).
      \begin{theorem}[{\cite{Ber:EES}}] \label{ML-S2.1} Let $H$ be a 
symmetric operator in a Hilbert space $\frak H$.
For the operator $H$ to be essentially self--adjoint in 
$\frak H$ it is 
necessary and sufficient that for some $b>0$ the
function $u=0$ is the only strong solution
of the Cauchy problems
      \begin{equation} \label{2.2}
\frac{du}{dt}(t)\pm (iH)^*u(t)=0,\qquad  t\in 
[0,b),\quad u(0)=0.
       \end{equation}
            \end{theorem}

We return to the operator $L=\ch^{-1}(J\frac{d}{dx}+B)$ from the previous 
Subsection \plref{sec32}. 

For a real number $\ga$ let $\ga_\pm(t)$ be the unique 
solution of the initial value problem
      \begin{equation}
         y'(t)= \pm c(y(t)),\quad y(0)=\ga.
    \label{ML-G2.1}
      \end{equation}
Here, $c$ is the function defined in \eqref{G2.3}.
Note that $c$ and $\frac 1c$ are locally integrable
and hence the first order equation \myref{ML-G2.1} with separated
variables has a unique solution. Moreover, if 
\begin{equation}
\pm\int_0^{\pm\infty}\frac{1}{c(x)}dx=\infty
\label{ML-G4.4}
\end{equation}
then the solutions $\ga_\pm(t)$ exist for all 
$t\in\R$.

\begin{prop}[Local energy estimate] 
Let $s_t$ be a strong solution of the equation
     \begin{equation}
   \frac{d}{dt}u-iL^* u=0
                      \label{ML-G2.0}
\end{equation}
defined for $|t|<\eps$. Moreover, assume that for some 
$\ga<\gb$ the
functions $\ga_\pm, \gb_\pm$ are defined for 
$|t|<\eps$.
Then the function
\[
    F_{\ga,\gb}(t):= \int_{\ga_+(t)}^{\gb_-(t)} 
s_t(x)^*\ch(x) s_t(x) dx
\]
is a decreasing function of $t$. 

In particular, if $\supp(s_0)\subset [\ga,\gb]$ then
$\supp(s_t)\subset [\ga_-(t),\gb_+(t)]$.
\end{prop}
\begin{proof}
     Differentiation by $t$ and integration by parts 
yields in
view of \eqref{ML-G3.26}
\begin{align}
   &    \frac{d}{dt}  \int_{\ga_+(t)}^{\gb_-(t)} 
s_t(x)^*\ch(x) s_t(x) dx
         \nonumber\\
   &= -c(\gb_-(t)) (s_t^*\ch s_t)(\gb_-(t))
      -c(\ga_+(t)) (s_t^*\ch s_t)(\ga_+(t))\nonumber \\
    &\quad -i  \int_{\ga_+(t)}^{\gb_-(t)} 
(L^*s_t)(x)^*\ch(x) s_t(x) - s_t(x)^* \ch(x) 
(L^*s_t)(x)dx\nonumber \\
   &= -c(\gb_-(t)) (s_t^*\ch s_t)(\gb_-(t))
      -c(\ga_+(t)) (s_t^*\ch s_t)(\ga_+(t))\nonumber \\
   &\quad  
       -i (s_t^*J^* s_t)(\gb_-(t))+i 
(s_t^*J^* s_t)(\ga_+(t))
       \label{ML-G2.2}
\end{align}
and by definition of $c$ this is $\le 0$ (cf. \eqref{ML-G3.5}). Note that
all terms in \myref{ML-G2.2} are real.

The last statement is clear.
\end{proof}

      \begin{prop}[Local existence]
 For each $f\in\cd_{\comp}(L^*)$ there exists an 
$\eps>0$ and a
unique strong solution $s_t, |t|<\eps,$ of the 
equation
\eqref{ML-G2.0} satisfying the initial condition $u(0)=f.$
Moreover $s_t\in\cd_{\comp}(L^*)$ for all $t$.
     \end{prop}
\begin{proof}
    Assume that $\supp(f)\subset [-N,N]$ for some 
$N>0$.
Choose a self--adjoint extension, $L_N$, of $L$ on the 
interval 
$[-2N,2N]$. This is possible since in view of Proposition
\plref{S1.9} the deficiency indices of $L$ on the finite
interval $[-2N,2N]$ are given
by $\cn_{\pm}(L)=N_\pm(L)=n$.

Next let $s_t(x):= e^{itL_N}f$ be the strong solution 
of the wave
equation for $L_N$. The local energy estimate above 
shows that
for $t$ small enough, $s_t$ has compact support in $[-
2N-\delta,2N+\delta]$
and hence can be extended by $0$ to a strong solution 
of
the wave equation for $L^*$.

The uniqueness follows immediately from the local energy estimate.
\end{proof}

Now we can give the

\begin{proof}[Second proof of Theorem \ref{ML-S1.3}]

If \eqref{ML-G1.6} is fulfilled then the previous 
result shows
that for each $f\in\cd_{\comp}(L^*)$ there exists a 
unique strong solution
$s_t, t\in\R$, of the Cauchy problem for the wave equation \eqref{ML-G2.0}
and $s_t\in\cd_{\comp}(L^*)$ for 
all $t$. Hence the result follows from Theorem 
\plref{ML-S2.1}.
      \end{proof}

      \section{Defect numbers and essential self--adjointness on the 
                  half-line}\label{sec5}

In this section we present some results on the 
square--integrable solutions of the system
   \begin{equation}\label{3.1}
J(x)y'(x)+B(x)y(x)=\gl {\ch}(x)y(x)
   \end{equation}
on the half lines $\mathbb R_{\pm}$. 
As in Section \ref{sec2} we associate with 
equation \eqref{3.1} the minimal symmetric linear 
relations $\csminpm$ 
and $\Sminpm$ in $\cl^2_{\ch}({\mathbb R}_{\pm},{\mathbb C}^n)$ 
and  
$L^2_{\ch}({\mathbb R}_{\pm },{\mathbb C}^n)$ 
respectively; $\cs_{\pm},S_\pm, \csmaxpm,\Smaxpm$ are
defined accordingly (cf. Def. \plref{ML-S2.0}).
As in Section \ref{sec2} we denote by
$\cn_{\pm}(\cs_{\pm}):=\dim \ce_{\pm i}(\cs_{\pm})$ the 
formal deficiency indices of the system \eqref{3.1}. 

If in addition ${\ch}(x)$ is invertible for almost all 
$x\in {\mathbb R}_{\pm}$ then $\cs_{\pm}$ is an 
operator. In this case the formal defect subspace 
$\ce_{\gl}(\cs_+)$ coincides with defect subspace 
$E_{\gl}(S_+)$ of the operator $S_+$.

We denote by $\kappa _+:=\kappa _+(iJ(0))$ and $\kappa _-
:=\kappa _-(iJ(0))$
respectively the numbers of positive and negative 
eigenvalues of the matrix $iJ(0)$. Since $\det J(x)\not= 0$ for $x\in I$
it is clear that $\kappa_{\pm}(iJ(x))$
does not depend on $x\in I, \kappa_{\pm}(iJ(0))= \kappa_{\pm}(iJ(x)).$
In what follows we will write sometimes $\kappa_{\pm}(iJ)$ instead of
$\kappa_{\pm}(iJ(0))$. 
Recall the well--known estimates (see {\cite[Theorem 9.11.1]{Atk:DCB})
   \begin{subequations}\label{3.2}
  \begin{align}
&\kappa _{\pm}\le \cn_{\pm}(\cs_+)\le n,\label{3.2a}\\  
&\kappa _{\pm}\le \cn_{\mp}(\cs_-)\le n,\label{3.2b}\\ 
&\cn_+(\cs_{\pm}) + \cn_-(\cs_{\pm})\ge n.\label{3.2c}
   \end{align}
       \end{subequations}
   \begin{remark}\label{rem2.6}
These inequalities have been established in 
\cite{Atk:DCB} by a 
generalization of the well--known Weyl analytic 
(circle--point) method. We note 
that in the case $\cn_+({\cs}_+)=\cn_-({\cs}_+),
\kappa_+=\kappa_-=n/2$ they 
follow easily from the results of Subsection \plref{sec25}:

For simplicity let us assume that the system $S_+$ is 
definite on $\R_+$.
Then by Proposition \plref{p3.1} (2),(3) 
$\dim({\cs}_+^*/{\cs}_+)\ge n$ 
and by \eqref{1.41} 
$\cn_+({\cs}_+) + \cn_-
({\cs}_+)=\dim({\cs}^*_+/{\cs}_+)\ge n.$ 
If $\cn_+({\cs_+})=\cn_-({\cs_+})$ then 
$\cn_{\pm }(S_+)\ge n/2.$ 
These inequalities
imply \eqref{3.2} if $\kappa _+(iJ)=\kappa _-(iJ)=n/2.$ 

However, we emphasize that we did not succeed to prove
the estimates 
\eqref{3.2} in full generality in the framework of 
extension theory.

Finally, note that, e.g., if $J^{-1}\ch$ is real
then $\cn_+(\cs_+)=\cn_-(\cs_+)$, 
cf. Proposition \plref{p3.4} below.
      \end{remark}
     \subsection{Minimal deficiency indices}\label{sec5.1}
Here we present a result on minimal possible deficiency 
(and formal deficiency) indices. It may be directly 
obtained by combining 
Theorem \plref{S2.2} and Proposition \plref{p3.2} below but in order to
demonstrate "formal" approach we present a simple independent proof.
     \begin{theorem}\label{th3.1} Let $c(x)$ be the function 
defined in \eqref{G2.3}. If 
  \begin{equation}\label{3.6}
\int ^{\infty }_0\frac{1}{c(x)}dx=\infty \qquad \Bigl(\text{resp. }\int ^0_{-
\infty }\frac{1}{c(x)}dx=\infty \Bigr)
   \end{equation}
then $\cn_{\pm}(\cs_+)=N_{\pm}(S_+)=\kappa _{\pm}$
(resp. $\cn_{\pm}(\cs_-)=N_{\pm}(S_-)=\kappa_{\mp}$).
     \end{theorem}
     \begin{proof} It suffices to prove the Theorem for the
linear relation $S_+$.
As noted in Remark \plref{ML-S3.3} it follows from \eqref{3.6} that 
$S_+$ is definite.
Therefore by Proposition \ref{p3.1} (1) 
$\cn_{\pm}(\cs_+)=N_{\pm}(S_+).$
Thus it suffices to prove the assertions for
$\cn_{\pm}(\cs_+).$

Let $y$ be a solution of \eqref{3.1} with $\lambda=\pm i.$ 
Let $(a_k)_{k\in\N}\subset\R_+$ be any sequence converging
to $\infty$. Then
integrating by parts and taking \eqref{G1.0} into account one gets

   \begin{equation}\label{3.7}
   \begin{split}
\|y\|^2_{\ch}&=\int_{\R_+}y(t)^*{\ch}(t)y(t)dt\\
     &=\lim_{k\to\infty}\Bigl(
         -\gl\int_0^{a_k}y^*(t)J(t)y'(t)dt-\gl\int_0^{a_k}y^*(t)B(t)y(t)dt\Bigr)\\
    &=\lim_{k\to\infty}\Bigl(\bigl[-\gl y^*(t)J(t)y(t)\bigr]\big|_0^{a_k} 
     \quad -\gl\int_0^{a_k}(J(t)y'(t))^*y(t)dt  \\
    &+ \gl\int_0^{a_k}y^*(t)J'(t)y(t)dt  - 
       \gl\int_0^{a_k}y^*(t)B(t)y(t)dt\Bigr)     \\
    &=\lim_{k\to\infty}\Bigl(\bigl[-\gl y^*(t)J(t)y(t)\bigr]\big|_0^{a_k}
        -\int_0^{a_k} y^*(t)\ch(t)y(t)dt\Bigr).
   \end{split}
   \end{equation}
Thus $\lim\limits_{k\to\infty} y(a_k)^*J(a_k)y(a_k)$ exists and
\begin{equation}
   2\|y\|_\ch^2=-\gl \lim_{k\to\infty}\bigl[y^*(t)J(t)y(t)\bigr]\big|_0^{a_k}.
   \label{ML-G5.5}
\end{equation}
On the other hand we find using \eqref{G2.3} 
      \begin{equation}
|y^*(t)J(t)y(t)|
\le c(t) \|\ch(t)^{1/2}y(t)\|^2.
\label{ML-G5.6}
     \end{equation}
We claim that there is a sequence
$(a_k)_{k\in\N}\subset\R_+$ such that $\lim\limits_{k\to\infty}
c(a_k)\|\ch(a_k)^{1/2}y(a_k)\|^2=0$. For if this were not the case
then we had an estimate $c(x)\|\ch(x)^{1/2}y(x)\|^2\ge \delta>0$ for
$x\ge x_0$. This would contradict \eqref{3.6} and
$\int_0^\infty \|\ch(x)^{1/2}y(x)\|^2dx=\|y\|_{\ch}^2<\infty$.

In view of \eqref{ML-G5.6} we have
  \begin{equation}\label{3.10}
\lim_{k \to \infty}\scalar{y(a_k)}{J(a_k) y(a_k)} = 0. 
   \end{equation}
Combining \eqref{ML-G5.5} and \eqref{3.10} 
one gets
   \begin{equation}\label{3.11}
2\|y\|^2_{\ch} = \scalar{\gl J(0)y(0)}{y(0)}_{\C^n}.
   \end{equation}
By the uniqueness theorem for first order differential equations
the map
$j: y(t)\to y(0)$
is an embedding of $\ce_\pm(\cs_+)$ into $\C^n$.
Moreover, the quadratic form $\scalar{iJ(0)\xi}{\xi}$ is positive
(resp. negative) on $j(\ce_+(\cs_+))$ (resp. $j(\ce_-(\cs_-))$).
Since $\kappa_\pm(iJ)$ is just the number of positive (resp. negative)
eigenvalues of the quadratic form $\scalar{iJ(0)\xi}{\xi}$ we 
obtain $\cn_\pm(\cs_\pm)\le \kappa_\pm$.
On the other hand we have in view of \eqref{3.2c}
\begin{equation}
 n\le \cn_+(\cs_+)+\cn_-(\cs_+)\le \kappa_++\kappa_-=n
\end{equation}
and thus equality holds. We emphasize that although we did not prove
\eqref{3.2} in full generality the relation 
\eqref{3.2c} was proved completely in Remark \plref{rem2.6}. 
      \end{proof}
      \begin{cor}\label{cor3.1}
Let $\gl_1(x)$ be the smallest eigenvalue of 
${\ch}(x)$. 
If for some $a\ge 0$
  \begin{equation}\label{3.13}
\int ^{\infty }_a\lambda _1(x)dx=\infty \qquad
(\int ^{-a}_{-\infty }\lambda _1(x)dx=\infty )
   \end{equation}
then $\cn_{\pm}(\cs_+)=N_{\pm}(S_+)=\kappa _{\pm}$
(resp. $\cn_{\pm}(\cs_-)=N_{\pm}(S_-)=\kappa_{\mp}$).
   \end{cor}
\begin{proof} This follows immediately from Corollary \plref{S2.3}.
   \end{proof}
     \begin{prop}\label{p3.2} Assume that the system 
\eqref{3.1} is definite on
${\mathbb R}_+$ and ${\mathbb R}_-$. Denote by 
$\cs, \cs_+, \cs_-, S, S_+$, and $S_-$ the symmetric linear relations 
associated to the equation \eqref{3.1} in 
$\cLHv{\R}, \cLHv{\R_{\pm}},  \LHv{\R}, \LHv{\R_\pm}$ respectively. Then
       \begin{subequations}\label{ML-G5.11}
\begin{align}
N_{\pm }(S)&=N_{\pm }(S_+)+N_{\pm }(S_-)-n,\label{ML-G5.11a}\\
\cn_{\pm}(\cs)&=\cn_{\pm}({\cs}_+)+\cn_{\pm}({\cs}_-)-n.\label{ML-G5.11b}
\end{align}
      \end{subequations}
            \end{prop}
\begin{proof} It follows from definiteness and Proposition \plref{p3.1} (1)
that \eqref{ML-G5.11a} and \eqref{ML-G5.11b} are equivalent. Hence it suffices
to prove one of them.

We put $\cs_0:={\cs}_+\oplus \cs_-$ 
and 
$S_0:=S_+\oplus S_-.$ 
By Proposition \plref{p3.1} (3) we have $f(0)=0$ for each 
$f\in \cd(\cs_0).$  
Moreover, Proposition \plref{p3.1} (2) implies that
for each $\xi\in {\C^n}$ there exists $\{f,g\}\in 
{\cs}$
with compact support such that $f(0)=\xi.$
Hence $\dim({\cs}/{\cs}_0)=n.$
In view of Proposition \plref{ML-S2.15} and Remark \plref{ML-S2.15a} the same
argument applies to $S_0$ and $S$. Hence $\dim(S/S_0)=n$.

On the other hand since $S$ is a closed symmetric 
extension 
of $S_0$ it follows from the second von Neumann formula
\eqref{1.43} with ${\widetilde S}$ and 
$S$ replaced by $S$ and $S_0$ respectively, 
that 
$N_{\pm}(S)=N_{\pm}(S_0)-\dim(S/S_0)=N_{\pm}(S_0)-n$.
Combining this formula with the obvious 
equalities $N_{\pm}(S_0)=N_{\pm}(S_+)+N_{\pm }(S_-)$ 
we obtain \eqref{ML-G5.11a} and thus also \eqref{ML-G5.11b}.
    \end{proof}
\begin{remark}\label{r3.1} 
1. The proof of the Proposition \plref{p3.2} is 
based essentially on the equality $\dim({\cs}/{\cs}_0)=n$ 
which is a consequence of Proposition \plref{p3.1}.
Note however that if ${\ch}$ is positive 
definite on $[-a,a]$ $(a>0)$ then this fact is obvious. 
Namely, $\{f_j,g_j\}^n_1$ forms a basis of 
${\cs}(\mod {\cs}_0)$ if 
$f_j\in \AC([-a,a],\C^n), \supp f_j\subset [-a,a], 
f_j(0)=\{\delta _{kj}\}^n_{k=1}$ and 
$g_j:=\chi {\ch}^{-1}(Jf'_j+Bf_j),\quad j\in 
\{1,...,n\}.$ Here, $\chi$ is a suitable cut--off function
with support in $[-a,a]$ and $\chi|\supp f_j=1$.

2. \eqref{ML-G5.11b} is due to Kogan and Rofe--Beketov 
\cite[Theorem 2.3]{KogRof:SIS}. Their proof is analytical in 
character and 
close to that given by Bennewitz \cite{Ben:GSS} for 
a similar formula for the scalar equation $Su=\gl Tu$, 
when one of the operators $S, T$ has a strictly 
positive Dirichlet 
integral on the solutions.

Our proof, being operator--theoretic in character, is 
rather simple and follows 
that of Glazman's result on ordinary differential 
equations on the line 
({\cite{AkhGla:TLO}}, {\cite{Naj:LDO}}).

3. Proposition \ref{p3.2} leads to a simple 
relation between Theorem \plref{S2.2} and 
Theorem \plref{th3.1}.
Indeed combining \eqref{ML-G5.11}, \eqref{3.2} and 
the obvious relation $\kappa_+ +\kappa_- = n$ we obtain 
the equivalences
\[\begin{split}
&N_+(S)=0 \Llr  N_+(S_{\pm })=\kappa _{\pm }, \\
&N_-(S)=0 \Llr  N_-(S_{\pm })=\kappa _{\mp}.
  \end{split}
\]
Thus Theorem \plref{S2.2} and Theorem \plref{th3.1} may be easily 
derived one from another.

4. \eqref{ML-G5.11} may be wrong for non--definite systems. For example
let
\begin{equation}
J=\begin{pmatrix}0&1\\-1&0\end{pmatrix},\quad
B=0,\quad\ch=\begin{pmatrix}1&0\\0&0\end{pmatrix}.
\end{equation}
One immediately checks that
\begin{equation}\begin{split}
     &N_+(S_\pm)=N_+(S)=N_-(S_\pm)=N_-(S)=0,\\
     &\cn_+(S_\pm)=\cn_+(S)=\cn_-(S_\pm)=\cn_-(S)=1.
		\end{split}
\end{equation}
Consequently, neither \eqref{ML-G5.11a} nor \eqref{ML-G5.11b} holds.
     
\comment{
 If the system \eqref{3.1} is not definite, then 
instead of \eqref{3.14}
one has 
   \begin{equation*}
n_{\pm}(S)=n_{\pm}(S_+)+n_{\pm}(S_-) + rank(S) - 
rank(S_+) - rank(S_-).
   \end{equation*}
This formula is implied by \eqref{3.16} and the 
obvious relations
$n_{\pm}(s)=N_{\pm}(S) + rank S-n.$}
   \end{remark}
\subsection{The case of singular Hamiltonian}\label{sec5.2}

Next we want to present a criterion for the deficiency indices to be minimal
on the half line $\R_+$ (for essential self--adjointness on the line $\R$)
in a case where the Hamiltonian is singular everywhere. 

We consider the type of first order systems
introduced in Example \plref{SDM-S1.2} and thereafter.
More precisely, we consider the first order system 
   \begin{equation}\label{ML-G3.11}
\tilde J f'+\tilde B f=\tilde \ch g,
  \end{equation}
where
   \begin{equation}\label{5.15A}
\tilde J=\begin{pmatrix} 0 & J^*\\ -J & 0\end{pmatrix}
,\quad
\tilde B=\begin{pmatrix} V & B \\ B^*-J' & -A \end{pmatrix},
\quad
\tilde\ch=\begin{pmatrix} \ch & 0\\0 & 0\end{pmatrix}.
      \end{equation}
$J,V,A,B,\ch$ are assumed to satisfy the same assumptions
as in \eqref{ML-G2.27}. In addition, $A$ is assumed to be non--negative.
Theorem  \plref{th3.1} does not apply to this situation
since $\tilde \ch$ is singular at every point. 
It is clear that $L_{\tilde \ch}^2(I)$ is canonically
isomorphic to $L^2_{\ch}(I)$.
We put (cf. \eqref{ML-G2.6})
$\tilde\cs_+=\cs_+(\tilde J,\tilde B,\tilde \ch)$.
For simplicity we will consider the interval $\R_+$ only.
For a function $f\in \cl_{\tilde \ch}^2(\R)$ we denote
by $f_1,f_2$ the first resp. last $n$ components.

We will use several times that if $\ch(x)$ and $A(x)$ are invertible
then we can estimate, for $\xi,\eta\in\C^n$,
   \begin{equation}
\begin{split}
   \big| \xi^* J(x)\eta|&\le\|A(x)^{1/2}\xi\|
   \|A(x)^{-1/2}J(x)\ch(x)^{-1/2}\ch(x)^{1/2}\eta\|\\
    &\le  \|A(x)^{-1/2}J(x)\ch(x)^{-1/2}\| \|A(x)^{1/2}\xi\| \|\ch(x)^{1/2}\eta\|.
\end{split}
\label{ML-G3.2.4}
   \end{equation}
Thus we put
   \begin{equation}
    c(x):=\begin{cases}  \max\bigl(1,\|A(x)^{-1/2}J(x)\ch(x)^{-1/2}\|\bigr),&
    \det(A(x)\ch(x))\not=0,\\
          \infty,&\textup{otherwise}.
          \end{cases}
\label{ML-G3.2.4a}
     \end{equation}
The self--adjointness criterion we are going to
present will depend also on $V$. We assume
that there exists an absolute continuous function $q\ge \delta>0$
on $\R$ such that
    \begin{equation}
    V\ge -q \ch.
\label{ML-G3.2.5}
   \end{equation}
      \begin{theorem}\label{th5.4} Let $A(x)$ be positive semi--definite for 
each $x\in \R_+$ and let $c(x)$ be the function defined in \eqref{ML-G3.2.4a}. 
Let $q\ge \delta >0$ be a function on $\R_+$ such that $V\ge -q{\ch}$  and
    \begin{equation}\label{5.43}
\int ^{\infty }_0\frac{1}{c(x) q^{1/2}(x)}dx=\infty \qquad
\bigl(\int ^0_{-\infty }\frac{1}{c(x) q^{1/2}(x)}dx=\infty\bigr).
    \end{equation}
Moreover, assume that one of the following two conditions is satisfied:
      \begin{thmenum}
\item $q$ is absolutely continuous and
\[\bigl|\frac{d}{dx}q^{-1/2}(x)\bigr|c(x)\le C_1\quad \text{ for }\quad x
\in \R_+;\]
\item $q(x)$ is non--decreasing (non--increasing).
\end{thmenum}
Then $\cn_{\pm }({\tilde{\cs}}_+)=N_\pm(\tilde S_+)=n$  
($\cn_{\pm }({\tilde {\cs}}_-)=N_\pm(\tilde S_-)=n$).
    \end{theorem}
    \begin{proof} 
The set
$\setdef{x\in\R}{\det(A(x)\ch(x))\not=0}$ has positive Lebesgue measure 
in view of  \eqref{ML-G3.2.4a}  and \eqref{5.43}.
Therefore by Proposition \plref{ML-S2.16}.
the system is definite. Hence it suffices to consider the formal
deficiency indices. 

1. Let $y$ be a solution of \eqref{ML-G3.11} with $\gl=\pm i$. We show that
$y_2\in \cl^2_{q^{-1}A}(\R_+)$. \eqref{ML-G3.11} reads
      \begin{equation}\label{5.44}
     \begin{split}
&J^*y'_2 + Vy_1 + By_2=\gl{\ch}y_1,     \\
&Jy'_1 - B^*y_1 + J'y_1 + Ay_2=0.
      \end{split}
   \end{equation}
It follows that 
      \begin{equation}\label{5.44'}
   \begin{split}
&\scalar{J^*y'_2}{y_1} + \scalar{Vy_1}{y_1}+\scalar{By_2}{y_1}=
\gl\scalar{{\ch}y_1}{y_1},\\
&\scalar{y_2}{Jy'_1} + \scalar{y_2}{J'y_1} - \scalar{By_2}{y_1} +
\scalar{y_2}{Ay_2}=0.
   \end{split}
  \end{equation}
Adding \eqref{5.44'}  and integrating from $0$ to $x$ one gets
       \begin{equation}\label{5.45}
     \begin{split}
F(x)^2 := &\int ^x_0q(t)^{-1}y^*_2(t)A(t)y_2(t)dt =- \int _0^x q(t)^{-1}
\scalar{y_2(t)}{J(t)y_1(t)}'dt\\
&-\int ^x_0q(t)^{-1}y^*_1(t)V(t)y_1(t)dt
+\gl\int ^x_0q(t)^{-1}y^*_1(t){\ch}(t)y_1(t)dt.
      \end{split}
          \end{equation}
We put 
$C_2=q(0)^{-1}|\Re\scalar{y_2(0)}{J(0)y_1(0)}|$
and recall (cf. \eqref{ML-G3.2.4}) that 
    \begin{equation}\label{5.46}
c^{-1}(x)|\scalar{y_2(x)}{J(x)y_1(x)}|\le \|{\ch}(x)^{1/2}y_1(x)\|\cdot 
\|A(x)^{1/2}y_2(x)\|. 
     \end{equation}
Using this and the inequality 
$|q(x)^{-3/2}q'(x)c(x)|\le C_1$ we obtain
    \begin{equation}\label{5.47}
    \begin{split}
\bigl|\int ^x_0\bigl(\frac{1}{q(t)}\bigr)'&\scalar{y_2(t)}{J(t)y_1(t)}dt|
\le \int ^x_0\frac{q'(t)}{q^2(t)}c(t)\|{\ch}^{1/2}(t)y_1(t)\|\cdot 
\|A^{1/2}(t)y_2(t)\|dt\\
\le &C_1\bigl(\int^x_0\|{\ch}^{1/2}(t)y_1(t)\|^2dt\bigr)^{1/2}\cdot 
\bigl(\int^x_0q(t)^{-1}\|A(t)^{1/2}y_2(t)\|^2dt\bigr)^{1/2}\\
\le &2^{-1}C_1^2\|y_1\|^2_{{\ch}} + 2^{-1}F^2(x).
   \end{split}
   \end{equation}
For brevity we assume in the sequel that $\delta =1$ that is $q(x)\ge 1.$
Now combining \eqref{5.46} and \eqref{5.47} and integrating by parts
we have
   \begin{equation}\label{5.48}
   \begin{split}
\bigl|\Re \int^x_0q(t)^{-1}\scalar{y_2(t)}{J(t)y_1(t)}'dt|\le &C_2 
+ c(x)\|{\ch}^{1/2}y_1(x)\|\cdot \|A^{1/2}y_2(x)\| \\
+ &2^{-1}C_1^2\|y_1\|^2_{{\ch}} + 2^{-1}F^2(x).
  \end{split}
  \end{equation}
Furthermore, the assumption $V\ge -q{\ch}$ yields
$-\int^x_0q^{-1}y_1^*Vy_1dt\le \int ^x_0y_1^*{\ch}y_1dt\le 
\|y_1\|^2_{{\ch}}.$ 
Thus setting
$C_3:=C_2+(2^{-1}C^2_1+1)\cdot \|y_1\|^2_{{\ch}}$ we infer from
\eqref{5.45} and \eqref{5.48} that
   \begin{equation}\label{5.49}
   \begin{split}
\int ^x_0\frac{F^2(t)dt}{2c(t)q(t)^{1/2}}&\le 
\int ^x_0\frac{C_3dt}{c(t)q(t)^{1/2}} +
\int ^x_0\frac{1}{q(t)^{1/2}}\|{\ch}^{1/2}y_1(t)\|\cdot \|A^{1/2}y_2(t)\|dt \\
&\le \int ^x_0\frac{C_3dt}{c(t)q(t)^{1/2}} + 
\|y_1\|_{{\ch}}\cdot F(x).
   \end{split}
  \end{equation}
We rewrite the latter inequality as 
    \begin{equation}\label{5.50}
G(x):=\int ^x_02^{-1}c^{-1}(t)q(t)^{-1/2}[F^2(t)-2C_3]dt
\le \|y_1\|_{{\ch}}F(x),
    \end{equation}
or as
   \begin{equation}\label{5.51}
G^2(x)\le \|y_1\|^2_{{\ch}}\bigl(2c(x)q^{1/2}(x)G'(x) + 2C_3\bigr).
   \end{equation}
We claim that $F^2(t)\le 2C_3$\ for $t\in {\mathbb R}_+.$ 
Assuming the contrary one finds $x_0$ such that $F(x_0)-2C_3=:\delta_1 >0$, 
hence $F(x)-2C_3\ge \delta_1$ for  $x\ge x_0$ since $F$ is non--decreasing.
Therefore in view of condition \eqref{5.43}  
$\lim\limits_{x\to \infty }G(x)=\infty.$

On the other hand choosing $a\in {\mathbb R}_+$ such that 
$G(a)\ge 2C_3^{1/2}\|y_1\|_{{\ch}}$, one derives from \eqref{5.51}
   \begin{equation*}
   \begin{split}
\frac{1}{2}\int^x_a\frac{dt}{c(t)q^{1/2}(t)}&\le 
\int ^x_a\frac{1}{c(t)q(t)^{1/2}}\bigl[1-2G(t)^{-2}
C_3\|y_1\|^2_{{\ch}}\bigr]dt  \\
&\le  \int ^x_a\frac{2G'(t)}{G^2(t)}dt = 2G(a)^{-1} - 2G(x)^{-1}
\le 2G(a)^{-1}.
    \end{split}
   \end{equation*}
This inequality contradicts the condition \eqref{5.43}.
Thus $q^{-1/2}y_2\in {\cl}^2_A({\mathbb R}_+)$ and 
$\|q^{-1/2}y_2\|^2_A\le 2C_3.$ 

2. Next we estimate using \eqref{ML-G3.2.4}
   \begin{equation}\label{5.53A}
   \begin{split}
|\scalar{y(x)}{\tilde J(x)y(x)}|& \le  2|\scalar{y_2(x)}{J(x)y_1(x)}|\\
     &\le c(x) \|A(x)^{1/2}y_2(x)\|\|\ch(x)^{1/2}y_1(x)\|\\
     &\le c(x) q(x)^{1/2} \|q(x)^{-1/2}A(x)^{1/2}y_2(x)\| \|\ch(x)y_1(x)\|.
   \end{split}
   \end{equation}
By 1. and Cauchy--Schwarz we know that $\|q(x)^{-1/2}A(x)^{1/2}y_2(x)\|
   \|\ch(x)y_1(x)\|$ is integrable. In view of the condition \eqref{5.43}
we infer exactly as in the proof of Theorem \ref{th3.1} that there
is a sequence $(a_k)_{k\in \N}\subset\R_+$ such that
$\lim\limits_{k\to\infty}|\scalar{y(a_k)}{\tilde J(a_k) y(a_k)}|=0$.
Also as in the proof of Theorem \ref{th3.1} one now completes
the proof, noting that $\kappa _{\pm }({\tilde J})=n.$

3. Now assume that condition (2) is satisfied. We reduce this case to the
previous one. For this purpose it suffices to construct an absolutely 
continuous function ${\tilde q}$ such that ${\tilde q}(x)\ge q(x)$ for 
$x\ge 0$ and ${\tilde q}$ satisfies both \eqref{5.43} and (1).

Since $c^{-1}(x)\le \|\ch(x)^{1/2}J^{-1}(x)A(x)^{1/2}\|\le 
\|\ch(x)^{1/2}\|\cdot
\|A(x)^{1/2}\|\cdot \|J^{-1}(x)\|$, one gets that 
$c^{-1} \in L^1_{\loc}({\mathbb R}_+).$ Therefore the function 
   \begin{equation}\label{5.27A}
t:= \gvf(x):=\int_0^xc(s)^{-1}ds
   \end{equation}
is absolutely continuous and monotone increasing for  $x>0.$ Denote by
$\psi$ the corresponding distribution function, 
$\psi(t):=mes\{x\in {\mathbb R}_+\ |\ \gvf(x)\le t\}.$ 

Next we put $q_1:=q\circ \psi$ and observe that $q_1$ is monotone 
increasing because so are $q$ and $\psi$. Besides it is clear that
    \begin{equation*}
\int _0^{\infty}q_1(t)^{-{1/2}}dt = \int_0^{\infty}q(x)^{-1/2}c(x)^{-1}dx =
\infty.
     \end{equation*}
Following F. S. Rofe-Beketov {\cite{Rof:SDO}} (see also {\cite{Shu:CQC}})
one puts  ${\tilde q}_1(n)=q_1(n+1)$ for $n\in {\mathbb Z}_+$ and then extends 
${\tilde q}^{-1/2}$ to the semi-axis  ${\mathbb R}_+$ by linear interpolation:
   \begin{equation*}
\tilde{q}_1(\gl n+(1-\gl)(n+1))^{-1/2}=
\gl\tilde{q}_1(n)^{-1/2} + (1-\gl)\tilde{q}_1(n+1)^{-1/2}, \quad \gl\in [0,1].
    \end{equation*}
It is clear that ${\tilde q}_1(x)\ge q_1(x)$ for $x\ge 0.$  
Moreover  ${\tilde q}_1^{-1/2}$ is globaly Lipschitz, 
$$
\left|
\frac d{dx}\tilde{q}_1(x)^{-1/2}\right|\le C_1:=q^{-1/2}(\psi(0))\quad 
\text{ and }\quad
\int_{\R_+}{\tilde q}_1(t)^{-1/2}dt=\infty.
$$
Finally, we put $\tilde{q}:=\tilde{q}_1\circ \gvf$ and check that $\tilde q$ 
has the desired properties. 

Indeed, $\tilde q(x)=\tilde q_1(\gvf(x))\ge 
q_1(\gvf(x))=q(\psi(\gvf(x)))\ge q(x)$, since $\psi(\gvf(x))\ge x$, and 
therefore $V\ge -{\tilde q}{\ch}.$
Further, ${\tilde q}^{-1/2}$ is absolutely continuous because so is
$\gvf$ and ${\tilde q}_1^{-1/2}$ is Lipschitz.   Now it follows from 
\eqref{5.27A} that 
$$
|({\tilde q}^{-1/2}(x))'|=
|({\tilde q}_1^{-1/2})'(\gvf(x))|
\cdot \gvf'(x)\le C_1 c(x)^{-1}
$$
and 
$$
\int_{\R_+}{\tilde q}(x)^{-1/2}c(x)^{-1}dx=
\int_{\R_+}{\tilde q}_1(t)^{-1/2}dt=\infty,
$$
which completes the proof. 
   \end{proof}
Combining Theorem \plref{th5.4} with Proposition \plref{p3.2} one arrives 
at the following self--adjointness criterion on the line. 
    \begin{theorem}\label{ML-S3.2.2}
Let $\tilde J,\tilde B,\tilde \ch$ be as in \eqref{5.15A}
with $A\ge 0$.
Let $q\ge \delta>0$ be a function on $\R$ such that $V\ge -q\ch$ and
      \begin{equation}\label{5.43'}
\displaystyle \pm \int_0^{\pm\infty} \frac{1}{c(x) q^{1/2}(x)} dx=
\infty.
     \end{equation}
Moreover, assume that one of the following two conditions is satisfied:
      \begin{thmenum}
\item $q$ is absolutely continuous and
\[\bigl|\frac{d}{dx}q^{-1/2}(x)\bigr|c(x)\le C_1\quad \text{ for }\quad 
x \in \R;\]
\item $q(x)$ is non--increasing on\ $\R_-$\ and  is non--decreasing on $\R_+$.
      \end{thmenum}
Then $\tilde \Smin=\Smin(\tilde J,\tilde B,\tilde \ch)$ is essentially 
self--adjoint.
   \end{theorem}
  \begin{remark}\label{r5.03} 
1. Let $V\ge 0$. In this case the proof of Theorem \plref{th5.4} essentially 
simplified and one easily gets that $y_2\in {\cl}^2_A({\mathbb R}_+)$ and
$\|y_2\|^2_A\le |\Re\scalar{J(0)y_2(0)}{y_1(0)}|.$

Moreover if $V\ge 0$ we may choose $q=1$. Then
Theorem \plref{th5.4}
holds under the only condition $\pm\int_0^{\pm\infty}\frac{1}{c(x)}dx=+\infty$.
   
2. The condition \eqref{5.43'}  is satisfied if $\|J(x)\|$ is bounded and
      \begin{equation}\label{5.43B}
\displaystyle \pm \int_0^{\pm\infty} \frac{\sqrt{\gl_1(A(x))\gl_1(\ch(x))}}
{q^{1/2}(x)} dx=\infty.
     \end{equation}
      \end{remark}


We apply Theorem \plref{th5.4} to the investigation of weighted matrix 
Sturm--Liouville (quasi--differential) equations with non--negative possibly 
singular (on some subsets of positive Lebesque measure) weight ${\ch }$
   \begin{equation}\label{5.53}
Py:=-\frac{d}{dx}\bigl(A(x)^{-1}\frac{dy}{dx}+Q(x)y\bigr) + 
Q^*(x)\frac{dy}{dx}+R(x)y=\gl \ch(x)y,
  \end{equation}
where we $A,Q,R,\ch$ satisfy the same assumptions as in Example 
\plref{SDM-S1.2}.

Denote by $\cn_{\pm }(P_+)$
the formal deficiency indices of the equation
\eqref{5.53} considered on the semiaxes ${\mathbb R}_+$, that is the number
of linearly independent solutions of \eqref{5.53} 
(with $\gl \in {\mathbb C}_{\pm})$  belonging to $\cl^2_{\ch}({\mathbb R}_+).$
By Proposition \plref{S1.11} the definition is correct, i.e. it does 
not depend on ${\pm}\gl \in {\mathbb C}_+$. 
    \begin{theorem}\label{th5.5} Let $P_+y=\gl {\ch}y$ be the equation of the 
form \eqref{5.53} with $A(x)$ being positive definite for 
$x\in {\mathbb R}_+,\ {\ch}\ge 0$ and $c(x)$ be defined by \eqref{ML-G3.2.4a}
with $J=iI.$
Suppose also that $V:=R-Q^*AQ \ge -q{\ch}$ where $q\ge \delta >0$ and 
    \begin{equation*}
\int ^{\infty }_0\frac{1}{c(x) q^{1/2}(x)}dx=\infty \qquad
(\int ^0_{-\infty }\frac{1}{c(x) q^{1/2}(x)}dx=\infty).
    \end{equation*}
Moreover, assume that one of the following two conditions is satisfied:
       \begin{thmenum}
\item $q^{-1/2}$ is absolutely continuous and

$\big|\frac{d}{dx}q^{-1/2}(x)\big|c(x)\le C_1\quad \text{ for }\quad
x\in {\mathbb R}_+$;
\item $q(x)$ is monotone increasing (monotone decreasing).
        \end{thmenum}
Then $\cn_{\pm }(P_+) = N_{\pm }(P_+)=n  \ 
\bigl(\cn_{\pm }(P_-) = N_{\pm }(P_-)=n\bigr).$
    \end{theorem}
   \begin{proof}

As elaborated in Example \plref{SDM-S1.2} the system \eqref{5.53} can be  
transformed into the first order system \eqref{ML-G3.11} 
$S(\tilde J, \tilde B, {\tilde \ch})$ with 
$\tilde J, \tilde B, {\tilde \ch}$ defined in \eqref{ML-G2.10}.
 
Namely, putting $u:=y$ and $v:=i(A^{-1}y'+Qy),$ one 
reduces the equation \eqref{5.53} to the  system 
     \begin{equation}\label{5.54}
\pmatrix 0&iI\\
iI&0
\endpmatrix
\binom{u}{v}'+
\pmatrix R-Q^*AQ&-iQ^*A\\
iAQ&-A
\endpmatrix
\binom{u}{v}=\lambda
\pmatrix{\ch}&0\\
0&0
\endpmatrix
\binom{u}{v}.
    \end{equation}
Since the corresponding linear relations are unitary equivalent,
we apply Theorem \plref{th5.4} and reach the conclusion.
      \end{proof}
    \begin{cor}\label{ML-S3.2.4}
Let $\Smin$ be the symmetric linear relation in $L^2_{\ch}(\R)$ induced
by the Sturm--Liouville type (quasi--differential) equation \eqref{SDM-G1.5}.
That is, $A,Q,R,\ch\in L^1_{\loc}(I,\Mat(n,\C))$, $A(x)$ is positive definite
for all $x\in \R$, and $\ch(x)\ge 0$. Let $c(x)$ be as defined in
\eqref{ML-G3.2.4a}.
Suppose that $V:=R-Q^*AQ\ge -q\ch$, where
$q\ge\delta>0$ and
     \begin{equation*}
\pm \int_0^{\pm\infty} \frac{1}{c(x) q^{1/2}(x)} dx=\infty.
    \end{equation*}
Let also one of the following two conditions is satisfied: 
       \begin{thmenum}
\item $q^{-1/2}$ is absolutely continuous and

$\big|\frac{d}{dx}q^{-1/2}(x)\big|c(x)\le C_1\quad \text{ for }\quad
x\in {\mathbb R}_+$;
\item $q(x)$ is monotone increasing on\ $\R_+$\  and is monotone decreasing
on $\R_-$.
        \end{thmenum}
Then $\Smin$ is essentially self--adjoint.
     \end{cor}
\begin{proof} This follows immediately from Theorem \plref{ML-S3.2.2} and
Example \plref{SDM-S1.2}.
     \end{proof}
  \begin{remark}\label{r5.04}
1. Another reduction of the equation \eqref{5.53} to the first order 
system has been used in {\cite{GohKre:TAV}} for the investigation of the 
asymptotic behavior of eigenvalues of boundary value problems for the 
equation \eqref{5.53}.

2. Theorem \plref{th5.5} generalizes some known results. 
Namely, for $Q=0$, $A={\ch}=I_n$ and real $R$ it has been obtained by
V. B. Lidskii {\cite{Lid:NSI}}. 
In turn for $n=1$  Lidskii's result coinsides with the well--known 
Titchmarsh-Sears theorem (see {\cite{BerShu:SE}}).

On the other hand, if $n=1, \quad Q=0,\ A=I_n$ and $R\ge 0$ the statement
of Theorem \plref{th5.5} has been established by M. G. Krein {\cite{Kre:TFO}}
(see also {\cite{KacKre:SFS}}). 
In Remark \plref{re5.40} below we will discuss also Krein's result 
for $R$ semibounded below $(R\ge -c\cdot I_n,\ c>0).$ 
   \end{remark}


   \subsection{Maximal deficiency indices}\label{sec3.3}
Here we investigate the opposite case of maximal 
deficiency indices.
  \begin{prop}\label{p3.3} Let $S_+=S_+(J,0,\ch)$ be a canonical system $(B=0)$
with  a Hamiltonian $\ch(x)=(h_{ij}(x))^{n}_{i,j=1}$ 
of positive type. If 
     \begin{equation} \label{3.17}
\int _0^{\infty }h_{jj}(x)dx<\infty,  \qquad \text{for } j=1,\ldots,k
      \end{equation}
then $\cn_{\pm}(\cs_+) = N_{\pm}(S_+)\ge \max \{{\kappa}_{\pm}, k\}$.
     \end{prop}
    \begin{proof} 
The condition $h_{jj}\in 
L^1(\mathbb R_+)$ is equivalent to the fact that the
constant vector $u_j:=\{\delta _{pj}\}^{n}_{p=1}$ is in 
${\cl}^2_{\ch}(\R_+, \C^n)$.
Thus $\ce_0(\cs_+) \supset \rmspan\setdef{u_j}{1\le j \le k}$ 
and $\dim \ce_0(\cs_+) \ge k.$ 
Since $\ch$ is of positive type the canonical system 
$S_+$ is definite. Therefore by Proposition \plref{p3.1} (1) we have
$\cn_{\pm}({\cs_+})=N_{\pm}(S_+)$ and
$\dim E_0(S_+)=\dim \ce_0({\cs}_+)\ge k.$ Now
Proposition \plref{p3.1} (7) implies the assertion.
    \end{proof}
   \begin{cor}\label{cor3.3} Let $S_+=S_+(J,B,\ch)$ be a 
definite system. Let $Y(x)=Y(x,0)$ be the fundamental
matrix solution of the equation \eqref{3.1} (cf. \eqref{G1.5})
and put ${\tilde {\ch}}(x):=Y^*(x){\ch}(x)Y(x)=
\bigl({\tilde h}_{ij}(x)\bigr)^n_{i,j=1}$. If the 
condition \eqref{3.17}
is satisfied with $h_{jj}$ replaced by ${\tilde 
h}_{jj}$, then
$\cn_{\pm }({\cs}_+)=N_{\pm }(S_+)\ge k$.
     \end{cor}
  \begin{proof} The 
gauge transformation $Y$ transforms the system into
a canonical one with Hamiltonian 
${\tilde {\ch}}$ and ${\tilde B}=0$ 
(see \eqref{ML1-G2.15}). 
A canonical system is definite if and only if the Hamiltonian
is of positive type. Hence
${\tilde{\ch}}$ is of positive type.
Since a gauge transformation preserves the deficiency indices 
we may apply Proposition \plref{p3.3} and reach the conclusion.
\end{proof}

\begin{theorem}\label{th3.2} Let $S_+=S_+(J,0,\ch)$ be a canonical system 
$(B=0)$
with a Hamiltonian $\ch$ of positive type on $\R_+$. 
For the equation \eqref{3.1} to have maximal formal 
deficiency indices 
$\cn_{\pm}(\cs_+)=n$ it is necessary and sufficient that 
   \begin{equation} \label{3.18}
\int^{\infty }_0\tr \ch(x)dx < \infty.
    \end{equation}
         \end{theorem}
    \begin{proof} 
\textit{Sufficiency.}
The inequality \eqref{3.18} is equivalent 
to \eqref{3.17} with $k=n,$ hence by Proposition 
\plref{p3.3} 
$\cn_{\pm}({\cs}_+)=n_{\pm}(S_+)\ge n$. 
On other hand $n\ge \cn_{\pm}({\cs}_+)$ and thus 
$\cn_{\pm}({\cs}_+)=N_{\pm}(S_+)=n$.

\textit{Necessity.}  Assume that $\cn_{\pm}({\cs}_+)=n.$ By 
Proposition \plref{p3.1} (1) also $N_{\pm}(S_+)=n$ 
and in particular $S_+$ admits self--adjoint extensions.
Fix one of them, say 
${\tilde S}_+={\tilde S}^*_+\supset 
S_+$. 

It follows from Proposition \plref{p3.1} (6) that 
there exists a linear relation
$\tilde \cs_+$ in $\cLHv{\R_+}$, satisfying
${\cs}_+\subset {\tilde {\cs}}_+\subset {\cs}^*_+$ 
and such that
$(\pi \oplus \pi ){\tilde {\cs}}_+={\tilde S}_+$. 
To calculate the
resolvent $({\tilde S}_+-\lambda )^{-1}$ we have to
find the solution
$\{{\tilde f},{\tilde g}\}\in {\tilde S}_+$ of the 
equation ${\tilde g}-\gl{\tilde f}={\tilde \psi }$ 
for an arbitrary
${\tilde \psi }\in \cLHv{\R_+}$, 
or what is 
the same, the solution $\{f,g\}\in {\tilde 
{\cs}}_+$ of the equation 
$Jf'-\lambda {\ch}f={\ch}\psi$  
with $f$ satisfying some (self--adjoint) boundary  
conditions at zero and at
infinity. It is well--known (see {\cite{Atk:DCB}}, 
{\cite{KogRof:SIS}}) that 
       \begin{equation}\label{3.19}
f(x,\lambda )=-\int ^{\infty }_0K(x,t,\lambda 
){\mathcal H}(t)\psi (t)dt=:
K_{\lambda }(\psi ),
        \end{equation}
where
   \begin{equation}\label{3.20}
K(x,t,\gl)=Y(x,\gl)[F(\gl)+1_{\R_+}(x-t)J^{-1}]Y(t,\ovl{\gl})^*.
   \end{equation}
Here $Y(x,\gl)$ is the fundamental $n\times n$ 
matrix solution of 
\eqref{3.1} (with $B=0$) satisfying the initial 
condition 
$Y(0,\gl)=I_n$ and $F(\gl)$ is some function.

It follows from \eqref{3.19} that $K_{\gl}(\psi)$ 
does not depend on the 
representative 
$\psi \in \cLHv{\R_+}$ of ${\tilde \psi}$. Thus 
$K_{\gl}$
is well defined on $\LHv{\R_+}$ and in view of 
\eqref{3.19}
         \begin{equation}\label{3.21}
{\tilde f}=({\tilde S}_+-\gl)^{-1}{\tilde \psi}=
\pi K_{\gl}(\psi) , \qquad
\gl \in \C_+\cup \C_-.
         \end{equation}
Combining \eqref{3.19}--\eqref{3.21} and $\cn_\pm(\cs_+)=n$
we see that the resolvent 
$({\tilde S}_+-\gl)^{-1}$ is a Hilbert-Schmidt 
operator for
$\gl \in {\mathbb C}\setminus {\mathbb R}.$ 
Consequently 
the spectrum 
$\sigma ({\tilde S}_+)$ is discrete.

Since $S_{\max,+}/S_+$ is finite--dimensional the existence
of a self--adjoint extension of $S_+$ with compact resolvent
implies that $S_+-\gl$ is a Fredholm relation of index 
$n$ for all $\gl\in\C$.
On the other hand by 
Proposition \plref{S1.8} (1) we have 
$\ker(S_+ - aI)=\{0\}$ for all $a\in \R$.
Therefore $\dim \ker(S_{\max,+}-aI)=n$. In particular 
$\dim \ker S_{\max,+}=n$ and by 
Proposition \plref{p3.1} (1) we obtain
$\dim\ce_0({\cs}_+)=\dim E_0(S_+)=n.$ 

But since the system is canonical we have
$E_0({\cs}_+)=\rmspan\{u_j\}^n_1$ with the constant vectors
$u_j=\{\delta_{pj}\}^n_{p=1}$. 
Thus $u_j\in \cLHv{\R_+}$ for $1\le j\le n.$ 
This is equivalent to 
$h_{jj}\in L^1(\R_+)$, $1\le j\le n,$ 
that is to the inequality 
     \eqref{3.18}.
         \end{proof}

To present the next result we recall the following definition.
  \begin{dfn}\label{def5.1} 
A symmetric system \eqref{3.1} is said to be quasi--regular 
if $\dim \ce_{\gl}({\cs}_+)=n$ for all $\gl \in \C,$ that is
$\cn_{\pm}({\cs}_+)=\dim \ce_a({\cs}_+)=n$ for all 
$a\in{\mathbb R}$.
   \end{dfn}
The following result is a refinement of Theorem \plref{th3.2}.
    \begin{theorem}\label{th3.3}
Under the conditions of Theorem \plref{th3.2} the system 
$S_+$
is quasi--regular on ${\mathbb R}_+$ if and only if  
$\int_{{\mathbb R}_+}\tr{\ch}(x)dx < \infty.$
\end{theorem}
\begin{proof} It is clear that $\cn_{\pm}({\cs}_+)=n$ 
if the system
$S_+$ is quasi--regular. Conversely, if 
$\cn_{\pm}({\cs}_+)=n$ 
then the relations $\dim \ce_a({\cs}_+)=n$ for 
$a\in {\mathbb R}$,
have been established in the proof of Theorem 
\plref{th3.2}.
\end{proof}
The next Corollary is derived from Theorem \plref{th3.3}
exactly as Corollary \plref{cor3.3} is derived from Proposition 
\plref{p3.3}.
    \begin{cor}\label{cor3.4} Let $S_+=S_+(J,B,\ch)$ be
definite on $\R_+$ and  
${\tilde{\ch}}$ be as in Corollary \plref{cor3.3}. 
Then for the 
system $S_+$ to be quasi--regular it is necessary and sufficient that
   \begin{equation}\label{5.80}
\int ^{\infty }_a\tr{\tilde {\ch}}(x)dx<\infty
    \end{equation}
    \end{cor}
      \begin{cor}\label{cor5.18A}
Let $S_+=S_+(J,B,\ch)$ be a definite system on $\R_+$ with constant
$J=J(0)$ and such that 
$\int_1^\infty x\|B(x)\|dx<\infty.$ Then for  the system $S_+$ 
to be quasi--regular it is necessary and sufficient that the 
condition \eqref{3.18} to be satisfied.
      \end{cor}
      \begin{proof}
It follows from the assumption  $\int_1^\infty x\|B(x)\|dx<\infty$
that there exists a fundamental $n\times n$ matrix 
solution $U(x)$ of the homogeneous equation 
$JU'(x)+B(x)U(x)=0$ satisfying 
     \begin{equation}\label{5.81}
U(x)=I_n + 0_n(1),\qquad x\to\infty
      \end{equation}
where $0_n(1)$ is $n\times n$ matrix function with entries $o(1).$
This fact is well known and can be easily checked (compare with the proof
of Proposition \plref{pr5.32A}).
By Corollary 5.17 $S_+$ is quasiregular iff 
$\int_0^\infty\tr(U^*(x)\ch(x)U(x))dx<\infty.$ 
In view of \eqref{5.81}
the last inequality is equivalent to the inequality \eqref{3.18}.
     \end{proof}
Another criterion for the formal deficiency indices 
$\cn_{\pm}$ to attain their
maximum values $n$ simultaneously (and thus a criterion 
for the system \eqref{3.1} to be quasi-regular) 
has been obtained in {\cite{KogRof:SIS}}:
   \begin{prop}\label{p3.5}{\cite[Theorem 
3.1]{KogRof:SIS}}
The system \eqref{3.1} is quasi-regular on $\R_+$ 
if and only if $\dim \ce_{\gl_0}({\cs}_+)=n$ 
for some $\gl_0\in \C$ and
   \begin{equation}\label{3.23}
\inf_{0\le t <{\infty}}\{\sgn(\Im(\gl_0))\int ^t_0
\tr(iJ(t)^{-1}{\ch}(t))dt\}>-\infty.
   \end{equation}
   \end{prop}
     \begin{remark}\label{r3.2} 
1. We emphasize that Theorem \plref{th3.2} as well as the other 
results of this subsection do not depend on $J$.

2. For Theorem \plref{th3.2} (as well as for 
Proposition \plref{p3.3}) to hold it is essential that ${\ch}$ 
is of positive type. Otherwise counterexamples are easy to find.
      \end{remark}
\subsection{Intermediate case}\label{sec5.4}
   \begin{dfn}\label{def5.2} Let $A$ be a linear 
relation in a Hilbert 
space $\frak H$ and let $j$ be an involution (that is an 
anti--linear bijective map) 
in $\frak H$. We will say that $A$ is invariant under 
$j$ if
$\{f,g\}\in A$ implies $\{jf,jg\}\in A$. 
    \end{dfn}
     \begin{lemma}\label{lem3.2} 
Suppose that the symmetric linear relation $A$ in
$\frak H$ is invariant under an involution $j$. Then 
$n_+(A)=n_-(A).$
      \end{lemma}
    \begin{proof} If $\{f,if\}\in \hat E_i(A)$ then 
$(j\oplus j)\{f,if\}=\{j f,-ij f\}\in A^*,$ hence 
$\{jf,-ij f\}\in \hat E_{-i}(A).$ Applying the same 
argument to $j^{-1}$ one sees that $j$ is an isomorphism
from $\hat E_{\pm}(A)$ onto $\hat E_{\mp}(A)$. 
    \end{proof}
   \begin{prop}\label{p3.4} Assume that $S_+=S_+(J,B,\ch)$
is definite on $\R_+$.  
If both $J^{-1}B$  and $J^{-1}{\ch}$ are real (that is 
have real entries) then
\begin{thmenum}
\item 
$N_+(S_+)=\cn_+({\cs}_+)=N_-(S_+)=\cn_-({\cs}_+);$
\item if $\dim \ce_{\lambda _0}({\cs}_+)=n$ for some
$\lambda_0\in \C$ then
    \begin{equation}\label{3.22}
N_{\pm }(S_+)=\cn_{\pm}({\cs}_+)=\dim 
\ce_a({\cs}_+)=\dim E_a(S_+)=n \quad 
\text{for any}\quad a\in {\mathbb R}..
    \end{equation}
\end{thmenum}
\end{prop}
   \begin{proof}
(1) ${\cs}_+$ is invariant under complex conjugation  
and therefore so is $S_+$. 
By Lemma \plref{lem3.2} $N_+(S_+)=N_-(S_-)$. 
The other equalities follow from 
Proposition \plref{p3.1} (1).

(2) If $\gl _0\in {\mathbb R}$ then the relations 
\eqref{3.22} are implied by Proposition \plref{p3.1} (7).
If $\lambda_0\in \C\setminus \R$ then by (1) 
$N_{\pm }(S_+)=\cn_{\pm }({\cs}_+)=n.$ The equality 
$\dim E_a(S_+)=n$ has been established in the 
proof of Theorem \plref{th3.2} (see also Theorem \plref{th3.3}).
      \end{proof}
  \begin{remark}\label{r3.3} 1. If $A$ is an operator 
then Definition \plref{def5.2} means that $A$ commutes 
with $j$. In this case Lemma \plref{lem3.2} is well--known.

2. The last three equalities in \eqref{3.22} meaning the 
quasi--regularity of the system
\eqref{3.1} have been established in 
{\cite[Theorem 9.11.2]{Atk:DCB}}
by an analytic method. A generalization of this result is 
contained in Proposition \plref{p3.5}.   
Note however that the condition 
\eqref{3.23}, meaning that the formal deficiency 
indices $\cn_{\pm}({\cs}_+)$
attain their maximum value simultaneously, does not 
imply the equality 
$\cn_+({\cs}_+)=\cn_-({\cs}_-)$ (see Example \plref{ex3.3} 
below).
\end{remark}
Now we are ready to present conditions for the canonical system 
\eqref{3.1}  
to have the formal deficiency indices $\cn_{\pm 
}({\cs}_+)=n-1.$
  \begin{prop}\label{p3.6}
Let $S_+$ be a canonical system on $\R_+$ with 
a Hamiltonian ${\ch}$ of positive type satisfying 
   \begin{equation}\label{3.24}
\int_{\R_+}h_{nn}(t)dt=\infty,\quad 
h_{jj}\in L^1(\R_+), \quad j=1,\ldots,n-1.
   \end{equation}
If in addition
   \begin{equation}\label{3.25}
\bigl|\int_{\R_+}\tr(iJ^{-1}{\ch}(t))dt\bigr|<\infty
   \end{equation} 
then $N_{\pm }(S_+)=\cn_{\pm }({\cs}_+)=n-1.$ 
   \end{prop}
\begin{proof}
Since $h_{jj}\in L^1(\R_+), j=1,\ldots,n-1$ 
then by
Proposition \plref{p3.3} $n-1\le N_{\pm }(S_+)\le n.$ 
Applying Theorem \plref{th3.2} we are, in view of condition
\eqref{3.24}, left with three 
possibilities:
  \begin{equation}\label{3.26}
N_{\pm }:=N_{\pm }(S_+)=n-1, \quad (N_+,N_-)=(n-1,n), 
\quad (N_+,N_-)=(n,n-1).
  \end{equation}
We rule out $N_-=n$ and $N_+=n.$  
The condition \eqref{3.25} yields \eqref{3.23}
with $\gl_0=-i$ and $\gl_0=i$. So if $N_-=n$ or $N_+=n$ then
by Proposition \plref{p3.5} the system \eqref{3.1} is 
quasi--regular, hence
$N_-=N_+=n.$ This contradicts \eqref{3.26}. Thus 
$N_{\pm }=n-1.$
      \end{proof}
         \begin{cor}\label{cor3.5} Let $S_+$ be a canonical system on $\R_+$ 
with a Hamiltonian ${\ch}$ of positive type such that 
$J^{-1}{\ch}$ is real. If the 
condition \eqref{3.24} is satisfied then 
$N_{\pm }(S_+)=\cn_{\pm }({\cs}_+)=n-1.$
    \end{cor}
      \begin{proof}
We show that the condition \eqref{3.25} is satisfied and apply Proposition
\plref{p3.6}. 
Since $J^{-1}{\ch}$ is
real so is $a:=\tr J^{-1}{\ch}.$ On the other hand 
$a=\tr(J^{-1}{\ch})=\tr({\ch}^{1/2}J^{-1}{\ch}^{1/2})
\in i\R$ since $J^{-1}$ is skew-adjoint. Thus $a=0$.
\end{proof}
 
In view of the importance of Hamiltonian systems we reformulate
Proposition \plref{p3.6} for such systems.
    \begin{cor}\label{cor3.6}
Let $n=2m, B=0, J=\begin{pmatrix} 0&I_m \\-I_m&0\end{pmatrix} $
and let ${\ch}=\begin{pmatrix} A&C\\C^*& D\end{pmatrix} $
be the block-matrix representation of a positive type 
Hamiltonian ${\ch}$
with respect to the decomposition $\C^n=\C^m\oplus \C^m.$
Suppose that the condition \eqref{3.24} holds and that
   \begin{equation}\label{3.27}
\bigl|\int _{\R_+}\tr(C_{I}(t))dt\bigr|<\infty, \qquad 
(C_{I}:=(C-C^*)/{2i}). 
   \end{equation}
Then $N_{\pm }(S_+)=\cn_{\pm}({\cs}_+)=n-1.$
   \end{cor}
    \begin{cor}\label{cor3.7} Let 
$J=\pmatrix 0&1 \\
-1&0
\endpmatrix$\   and ${\ch}=
\pmatrix a&b \\
{\overline b}&c
\endpmatrix$ 
be a $2\times 2$ Hamiltonian satisfying
$\bigl|\int_{\mathbb R_+}b_I(t)dt\bigr| <\infty.$ 
Moreover assume that the system $S_+$ is definite and
$\int_1^{\infty}x\|B(x)\|dx < \infty.$
Consider the symmetric extensions of $S_+$ defined by
     \begin{equation}\label{3.28}
{\widetilde S}_i:=\bigsetdef{\{{\tilde f}, {\tilde 
g}\}\in S_{\max,+}}{f=\col(f_1,f_2)\in \cd(\cs_{\max,+}),   
f_i(0)=0},\quad (i=1,2).
      \end{equation}
Then
\begin{thmenum}
\item $N_{\pm }({\tilde S}_i)=\cn_{\pm }({\tilde 
{\cs}}_i)=1$ if and only if
$\int_{\R_+}\tr{\ch}(x)dx < \infty$.

\item ${\tilde S}_i$ is self--adjoint, 
i.e. $\cn_{\pm}({\tilde {\cs}}_i)=N_{\pm}(\tilde S_i)=0$
if and only if $\int_{\R_+}\tr{\ch}(x)dx = \infty.$
\end{thmenum}
   \end{cor}
    \begin{proof} Since the system $S_+$ is definite 
then $\cn_{\pm}({\tilde {\cs}}_i)=N_{\pm}(\tilde S_i) $.
It follows from Proposition \plref{p3.1} (3) and \eqref{3.28} that 
$\dim({\widetilde \cs}_i/\cs_+)=1.$ By Proposition \plref{p3.1} (6)
we have $\dim({\widetilde S}_i/S_+)=1$, too.
Hence $N_\pm(\tilde S_i)=N_\pm(S_+)-1\le 1$.

1. Let $\int_{\R_+}\ch(x)dx<\infty$. Then
by Corollary \plref{cor5.18A}
$N_{\pm }(S_+)=2.$ and thus $N_{\pm}({\tilde S}_i)=1.$

2. Conversely, assume that $N_+(\tilde S_i)=1$ or $N_-(\tilde S_i)=1$.
Then $N_+(S_+)=2$ or
$N_-(S_+)=2$. As in the proof of Proposition \plref{p3.6} one
now concludes that the system is quasi--regular and hence
$N_+(S_+)=N_-(S_+)=2$.
\end{proof}
         \begin{remark}\label{r3.4}
Corollary \ref{cor3.7} slightly improves a result due to Kac--Krein 
\cite{KacKre:SFS} and
coincides with it if $B=0$ and $b={\overline b}$, that is $b_I=0.$ 
Our Theorem \ref{th3.2} has been inspired by this result.

Note also that the equalities $N_{\pm}(S_+)=1$ for $2\times 2$ definite 
systems with real trace-normed Hamiltonian ($\tr \ch(x)=1$ for $x\in \R_+)$ 
has been established by de Branges \cite{Bra1:SHS}. Another proof of the 
de Branges result has been proposed in the recent publication 
\cite{HSW:BVP}. These authors have also established an interesting 
inequality: 
    \begin{equation*}
{\scalar{f(x)-f(y)}{f(x)-f(y)}}_{{\mathbb C}^2}
\le \sqrt{6}|\gl|\sqrt{|x-y|}\cdot\|f\|_{\ch}\quad \text{for}
\quad f\in E_\gl(S_+).
      \end{equation*}
  \end{remark}

Now we present some examples clarifying the sharpness 
of the conditions 
\eqref{3.24} and \eqref{3.25} in Proposition 
\plref{p3.6}.
    \begin{example}\label{ex3.1} Let $J=\diag(i,-i),$
$\ch=\diag(h_{11},h_{22})$ where $h_{jj}(x)>0$ for  
$x\in \R_+$. 
If $h_{11}\not \in  L^1(\R_+)$ and $h_{22}\in 
L^1(\R_+)$ then 
the condition \eqref{3.24} holds but the 
condition \eqref{3.25} fails. 
It is easily seen  that $N_+(S_+)=\cn_+({\cs}_+)=1$ and 
$N_-(S_+)=\cn_-({\cs}_+)=2.$
If conversely $h_{11}\in L^1(\R_+)$ and 
$h_{22}\not \in L^1(\R)$
then $N_+=\cn_+=2$ and $N_-=\cn_-=1.$

This example shows that generally speaking Corollary 
\plref{cor3.7} does not occur if the condition 
\eqref{3.25} fails.
   \end{example}
     \begin{example}\label{ex3.2}
1. Let $J$ and ${\ch}$ be as in the previous 
example. Suppose that
$h_{11}(x)\ge h_{22}(x)>0$ for $x\in\R_+$,
$h_{22}\not \in L^1(\R_+)$ and 
$h_{11}-h_{22}\in L^1(\R_+).$ Then
$N_{\pm }=\cn_{\pm }=1$ though the condition 
\eqref{3.24} fails and the
condition \eqref{3.25} holds.

This example shows that the condition \eqref{3.24} is 
not necessary for the
relations $\cn_{\pm}({\cs}_+)=n-1$ to be valid.

2. If $h_{11}-h_{22}\not \in L^1(\R_+)$ (say 
$h_{11}=2(1+x)^{-1}, h_{22}=(1+x)^{-1}$)
then again $N_{\pm }=\cn_{\pm }=n-1=1,$ but neither 
condition \eqref{3.24} nor condition \eqref{3.25} hold.
    \end{example}
\begin{example}\label{ex3.3} We put $J=J_1\oplus 
J_1\oplus J_1$, where $J_1=\diag(i,-i)$, 
${\ch}=\diag(h_{11},...,h_{66}),$ and
$h_{11}=h_{33}=2^{-1}h_{66}\not \in L^1(\R_+)$ and 
$h_{22}=h_{44}=2^{-1}h_{55}\in L^1(\R_+)$. 
It follows from Example \plref{ex3.1} that 
$\cn_+({\cs}_+)=4$ and $\cn_-({\cs}_+)=5.$

On the other hand $\tr(J^{-1}{\ch})=0$ 
and hence the condition
\eqref{3.23} holds. 
This example shows that the condition 
\eqref{3.23} is not sufficient for the system \eqref{3.1} 
to have equal formal deficiency indices.
     \end{example}
\begin{example}\label{ML-S5.33} 
We put in Corollary \ref{cor3.7}  
$b(x)=0,$ $a(x)=(1+x)^{-4}, \quad c(x)=1.$  Then by Corollary
\plref{cor3.7} the operator $S_i$ is self--adjoint. 

On the other hand the eigenvalues of $\ch^{1/2}(x)J\ch^{1/2}(x)$
are $\pm i(1+x)^{-2}$. Hence we infer
that $c^{-1}(x)=(1+x)^{-2}\in L^1({\mathbb R}_+).$ 

This example shows that the conditions of Theorem \plref{th3.1} 
(Theorem \plref{S2.2}) are not necessary for $S_+$ to 
have minimal deficiency indices (to be self--adjoint).

Moreover, this example (as well as Example \plref{ex3.6}) shows that 
$S_i=S_{\max,i}$ though $\cd(\Smax)$ is not contained in 
$\AC_0({\mathbb R}_+,{\mathbb C}^2)$ (cf. Remark \plref{rem2.2'}). 
Indeed, put
     \begin{equation*}
f=\col((1+x)^{1/4},0), \quad  g=\col(0,-\frac14 (1+x)^{-3/4})\in 
\cLHv{\R_+}.
      \end{equation*}
Then $\{f,g\}\in \csmax$ and
${\scalar{f(x)}{f(x)}}_{{\mathbb C}^2} = \sqrt{1+x}\to 
\infty$ 
as $x\to \infty.$
        \end{example}
   \subsection{Two-terms Sturm-Liouville equation}\label{sec5.5} 
Let  us consider the equation \eqref{SDM-G1.5} with $Q=R=0,$ that is
   \begin{equation}\label{5.35}
Py:=-\frac{d}{dx}\bigl(A^{-1}\frac{dy}{dx}\bigr)=\gl \ch y.
   \end{equation}
      \begin{prop}\label{ML-S5.16}
Let $A(x)$ be positive definite for all $x\in \R_+$ and
$\ch(x)\ge 0$ and let $\ch(x)$ be a nonsingular on a subset of positive
Lebesgue measure.
Then for the equation \eqref{5.35} to have maximal formal 
deficiency indices $\cn_{\pm }(P_+)=2n$ (as well as to be quasiregular),
it is necessary and sufficient that 
  \begin{equation}\label{5.36}
\int^{\infty}_0\tr\bigl({\widetilde A}(x)\ch(x){\widetilde A}(x)\bigr)dx<\infty
\quad\text{and}\quad\int^{\infty }_0\tr(\ch(x))dx<\infty,
  \end{equation}
where ${\widetilde A}(x):=\int^x_0 A(t)dt.$ 

If $A$ is uniformly definite on 
$\R_+,$ that is $A(x)\ge \eps \cdot I$ $(x\in \R_+)$ with some 
$\eps >0$ then the second condition in \eqref{5.36} is obsolete. 
         \end{prop}
    \begin{proof}
As explained
in Example \plref{SDM-S1.2} the system $P$ is unitarily equivalent
to a first order system $S(\tilde J,\tilde B, \tilde \ch)$, with
$\tilde J, \tilde B\,\tilde \ch$ defined in \eqref{ML-G2.10}.
By Proposition \plref{ML-S2.16}    
the system $S(\tilde J,\tilde B, \tilde \ch)$ is definite.
Then the gauge transformation 
$Y=
\pmatrix I&-i{\tilde A}\\
0&I
\endpmatrix
$
transforms the system $S({\tilde J},{\tilde B},{\tilde \ch})$
into a canonical (and definite) one $S({\tilde J}, 0, {\tilde {\ch}_1})$ 
with ${\tilde J}$ and ${\tilde {\ch}_1}$ defined by
    \begin{equation}\label{5.37}
\tilde J=\pmatrix 
0&iI\\
iI&0
\endpmatrix\qquad \mbox{ and }\qquad 
{\tilde \ch}_1=Y^*{\tilde H}Y=
\pmatrix 
\ch&-i\ch \tilde A\\
i\tilde A\ch&\tilde A\ch\tilde A\endpmatrix.
    \end{equation}

Since the Hamiltonian ${\tilde \ch _1}$ is of positive type the first 
assertion follows from Theorem \plref{th3.2}.

To prove the second assertion we put $\ch _1:={\tilde A}\ch {\tilde A}$  
and $\ch _2:=\ch ^{1/2}{\tilde A}^2\ch ^{1/2}.$ 
Since $A(x)>\eps \cdot I$ one gets $\ch _2(x)\ge (\eps x)^2\ch (x).$ 
Using this and the equality $\tr \ch_1(x)= \tr \ch_2(x)$ 
we get
     \begin{equation*}
\int ^{\infty }_1\tr\ch _1(x)dx=\int ^{\infty }_1\tr\ch _2(x)dx\ge 
{\eps}^2 \int^{\infty }_1x^2\tr\ch (x)dx\ge 
{\eps}^2\int^{\infty }_1\tr\ch (x)dx.
      \end{equation*}
This proves the last statement. 
      \end{proof}

Similarly, starting with Proposition \plref{p3.3} and taking 
\eqref{3.2a} into account one arrives at the following 
   \begin{prop}\label{pr5.33}
Assume that the conditions of Proposition \plref{ML-S5.16} are fulfilled and 
$\ch=:(h_{ij})_{i,j=1}^n$ and 
$\tilde A\ch \tilde A=({\tilde h_{i,j}})_{i,j=1}^n.$ 
If 
   \begin{equation}
\int_0^\infty h_{jj}(x)dx<\infty,\ \ j\in \{1,\ldots,k_1\}\quad\text{and}\quad
\int_0^\infty {\tilde h}_{ii}(x)dx<\infty,\ \ i\in \{1,...,k_2\}
     \end{equation}
then $\cn_{\pm}(P_+)\ge \max \{n, k_1+k_2\}.$ 
      \end{prop}
    \begin{cor}\label{ML-S5.17}
Let $0<c_1\le A(x)\le c_2$ for $x\in \R_+$  and let $\ch(x)$ be positive 
definite on a subset of positive Lebesgue measure.
Then for the equation 
\eqref{5.35} to have maximal formal deficiency indices 
$\cn_{\pm }(P_+)=2n$ it is necessary and sufficient that
   \begin{equation*}
\int^{\infty }_0 x^2 \tr\ch(x)dx<\infty.
   \end{equation*}
    \end{cor}

Next we slightly generalize Proposition \plref{ML-S5.16}. Consider the 
matrix equation \eqref{SDM-G1.5} with $Q=0,$ that is 
   \begin{equation}\label{5.60}
Py:=-\frac{d}{dx}\bigl(A^{-1}\frac{dy}{dx}\bigr) + R(x)y = \gl \ch y.
   \end{equation}
     \begin{prop}\label{pr5.32A}  Assume that $\ch(x)$ is positive definite
on a subset of positive Lebesgue measure and
    \begin{equation}\label{5.60A}
\int_1^{\infty }\|{\tilde A}(x)\|\cdot\|R(x)\|dx<\infty \quad \text{and}\quad 
\lim_{x\to \infty} A(x)\int_x^{\infty}R(t)dt = 0.
  \end{equation}
Then for the equation \eqref{5.60}  to have maximal formal deficiency indices 
$\cn_{\pm }(P_+)=2n$ (as well as to be quasiregular) it is necessary and 
sufficient that the conditions \eqref{5.36} be satisfied.
    \end{prop}
  \begin{proof} At first we prove that the homogeneous equation 
\eqref{5.60} (with $\gl =0$) has two $n\times n$ matrix solutions $U$ and 
$V$ satisfying:
      \begin{equation}\label{5.61}
U(x)= I_n + 0_n(1), \qquad  U'(x)= 0_n(1), \quad  x\to \infty, 
   \end{equation}
     \begin{equation}\label{5.62}
V(x)={\tilde A}(x)\cdot \bigl(I_n + 0_n(1)\bigr), \quad 
V'(x)=A(x)\cdot \bigl(I_n + 0_n(1)\bigr), \quad  x\to \infty 
      \end{equation}
where as before $0_n(1)$ stands for the $n\times n$ matrix function
with entries $o(1)$ as $x\to \infty$.
Indeed it is clear that each solution $U$ of the integral equation
       \begin{equation}\label{5.63}
U(x)= I_n + \int_x^{\infty}A(t)dt\int_t^{\infty }R(s)U(s)ds
      \end{equation}
is also a solution of the equation \eqref{5.60} with $\gl =0.$ 
Choose $N$ such that
      \begin{equation}\label{5.64}
\int_N^{\infty}\|{\tilde  A}(s)\|\cdot\|R(s)\|ds < 1/2.
      \end{equation}
Further, setting  $U_0(x)=I_n$  and 
      \begin{equation*}
U_n(x)=\int_x^{\infty }A(t)dt\int_t^{\infty}R(s)U_{n-1}(s)ds 
= \int_x^{\infty}[{\tilde A}(s)-{\tilde A}(x)]R(s)U_{n-1}(s)ds  \quad (n\ge 1)
      \end{equation*}
and using \eqref{5.64} and the inequality
     \begin{equation*}
\|({\tilde A}(s)-{\tilde A}(x))R(s)U_{n-1}(s)\| \le 
\|{\tilde A}(s)\|\cdot \|R(s)\|\cdot \|U_{n-1}(s)\|, \quad  s>x,
   \end{equation*}
one easily proves by induction that 
$\|U_n(x)\|\le 1/2^n$ for $n\ge 1.$ Hence the series
$\sum^{\infty }_{n=1}U_n(x)$ converges uniformly for $x\ge N$ and 
$\|\sum^{\infty }_{n=1}U_n(x)\|\le 1.$ Moreover, the matrix function
$U(x):=I_n + \sum_{n\ge 1}U_n(x)$
defines the unique solution of the equation \eqref{5.63} (for $x\ge N$) and 
satisfies the inequality $\|U(x)\|\le 2$\ for\ $x\ge N.$

Using this estimate one obtains from \eqref{5.63} that 
$U(x)-I_n=0_n(1)\ \text{as}\  x\to \infty .$ Differentiating \eqref{5.63} 
and applying \eqref{5.60A} and the above estimate $\|U(x)\|\le 2$ one 
derives the second relation $U'(x)=0_n(1)\ \text{as}\  x\to \infty .$ 

Thus the existence of the solution $U$ satisfying \eqref{5.61} is proved.

To prove the existence of the solution $V$ satisfying \eqref{5.62}
we recall (see \cite{{Har:ODE}}, part XI) that for each $n\times n$
matrix solution of the equation \eqref{5.60} (with $\gl =0)$ the matrix 
function
   \begin{equation}\label{5.64A}
K:=U^*(x)A^{-1}(x)U'(x)- \bigl(A^{-1}(x)U'(x)\bigr)^*U(x)
   \end{equation}
is constant. Turn $x\ \text{to}\ +\infty$ and taking \eqref{5.61}
into account one gets that $K=0$. This means that $U$ is a self-adjoint
solution (in the sence of {\cite{Har:ODE}}, part XI) of the 
homogeneous equation \eqref{5.60} (with $\gl =0$).

Using \eqref{5.64A} (with $K=0$) it is easy to check (and it is known 
(see {\cite{Har:ODE}}), that the classical Liouville formula remains 
valid for the matrix case, that is
   \begin{equation}\label{5.64B}
V(x):=U(x)\int_0^xU^{-1}(t)A(t)(U^{-1})^*(t)dt
    \end{equation}
is also a $n\times n$ matrix solution of the equation \eqref{5.60} 
(with $\gl =0$). 

The relations \eqref{5.62} are implied now by \eqref{5.61} and \eqref{5.64B}.

Further, following the proof of Proposition \plref{ML-S5.16} one transforms 
the system $P$ to a first order system $S(\tilde J,\tilde B,{\tilde \ch})$ 
with $\tilde J,\tilde B,{\tilde \ch}$ defined in \eqref{ML-G2.10}. 
Then the gauge transformation 
$Y=
\pmatrix U&-iV\\
iA^{-1}U'&A^{-1}V'
\endpmatrix
$
transforms the system $S(\tilde J, \tilde B, {\tilde \ch})$
into a canonical system $S(\tilde J_1, 0, {\tilde \ch}_1)$
with 
   \begin{equation}\label{5.69}
\tilde J_1 = Y^*(0){\tilde J}Y(0)\quad \text{and} \quad 
{\tilde \ch}_1=Y^*{\tilde \ch}Y=
\pmatrix U^*{\ch}U&-iU^*{\ch}V\\
iV^*{\ch}U&V^*{\ch }V.
\endpmatrix
     \end{equation}
We note that generally speaking ${\tilde {J_1}} \not = {\tilde J}$ 
since $Y(0)\not =I.$

By Theorem \plref{th3.2} $N_{\pm }(P_+)=2n$\  
iff\  $\int_1^{\infty}\tr(U^*{\ch}U + V^*{\ch}V)dx<\infty.$ 
By Theorem \plref{th3.3} this inequality is also equivalent to the property
of the system $P_+$ to be quasiregular. In view of \eqref{5.61} 
and \eqref{5.62} this inequality is equivalent to \eqref{5.36}). 
      \end{proof}
   \begin{remark}\label{re5.37A}
1. If $\|A(x)\|$ is bounded $(\|A(x)\| \le C),$ then both conditions 
\eqref{5.60A} are implied by the condition 
$\int_1^{\infty }x\|R(x)\|dx < \infty.$ 

2. In the scalar case $(n=1)$ the second condition in \eqref{5.60A}
may be omitted.
   \end{remark}

Next we consider the equation \eqref{5.60} with $A=I.$ For this case
we complement Proposition \plref{pr5.32A}. 
    \begin{prop}\label{pr5.33A}
Let $A=I$ and let $\ch(x)$ be nonsingular on a subset of positive Lebesgue
measure.  Assume also that $R(x)=k^2\cdot I_n + R_1(x)$ where 
$\int_0^{\infty} \|R_1(x)\|dx < \infty .$ 
Then for the equation \eqref{5.60} to have
maximal formal deficiency indices $\cn_{\pm }(P_+)=2n$ 
(as well as to be quasiregular) it is necessary and sufficient that:
    \begin{equation}
      \begin{split}
&i)\qquad\qquad \int ^{\infty }_0\tr \ch(x)dx < \infty  \quad \text{if} \quad 
k=is\in i\R\ \ (k\not =0); \\
&ii)\qquad \qquad\int ^{\infty }_0e^{2kx}\tr \ch(x)dx < \infty  
\qquad \text{if}\ \ k>0.
     \end{split}
    \end{equation}
 \end{prop}
  \begin{proof}
i) If $||R_1||\in L^1(\R_+)$ then, as it is well known, there exist two
$n\times n$ matrix solutions $U$ and $V$ of the homogeneous equation
$-y''-s^2y + R_1(x)y=0$ satisfying
\begin{equation}
   \begin{split}\label{5.70}
U(x)&=\cos sx \cdot \bigl(I_n + 0_n(1)\bigr),\ \ U'(x)=-s \sin sx \cdot
\bigl(I_n + 0_n(1)\bigr), \ \ x \to \infty,  \\
V(x)&=\frac{\sin sx}{s} \cdot \bigl(I_n + 0_n(1)\bigr),\  \ 
V'(x)=\cos sx \cdot \bigl(I_n + 0_n(1)\bigr), \quad x \to \infty.
 \end{split}
        \end{equation}
Following the proof of Proposition \plref{pr5.32A} and using the gauge 
transformation
$Y=
\pmatrix U&-iV\\
iU'&V'
\endpmatrix
$
we reduce the equation \eqref{5.60} to a canonical system 
$S({\tilde J}_1, 0, {\tilde \ch }_1)$ with ${\tilde J}_1$ and 
${\tilde \ch}_1$ defined in \eqref{5.69}. 
In view of  \eqref{5.70} the inequality
$\int_0^{\infty}\tr(U^*\ch U + V^*\ch V)dx <  \infty$\ takes place iff\
$\int^{\infty}_0\tr\ch (x)dx < \infty .$
It remains to apply  Theorem \plref{th3.2}.

ii) Now the homogeneous equation $-y''+k^2y+R_1(x)y=0$ has two $n\times n$
matrix solutions satisfying
\begin{equation}
  \begin{split}
U(x)&=\cosh kx\cdot \bigl(I_n + 0_n(1)\bigr), \quad U'(x)=
k\cdot \sinh kx\cdot \bigl(I_n + 0_n(1)\bigr), \ \ x \to \infty, \\
V(x)&=k^{-1}\sinh kx\cdot \bigl(I_n + 0_n(1)\bigr),\quad V'(x)=
\cosh kx\cdot \bigl(I_n + 0_n(1)\bigr),\ \ x \to \infty. 
  \end{split}
           \end{equation}
Starting with these solutions one completes the proof in just the same way as
in the case i).
  \end{proof}
Next we present few results on intermediate formal deficiency indices
$\cn_{\pm}(P_+).$

   \begin{prop}\label{pr5.38} 
Let $\ch=:(h_{ij})_{i,j=1}^n$ and 
$\tilde A\ch\tilde A=:({\tilde h}_{ij})_{i,j=1}^n.$ 
Assume also that all the functions $\{h_{jj},{\tilde h}_{jj}\}_1^n$ 
but one belong to the space $L^1(\R_+).$ Then the formal deficiency 
indices of the equation \eqref{5.35} are $\cn_{\pm}(P_+)=2n-1.$ 
  \end{prop}
    \begin{proof}
As in the proof of Proposition \plref{ML-S5.16} we transform the equation 
\eqref{5.35} to a canonical system $S(\tilde J,0,{\tilde\ch}_1)$ with 
$\tilde J$ and ${\tilde \ch}_1$ defined in \eqref{5.37}. 
One checks that  $\tr({\tilde J}^{-1}{\tilde \ch}_1(t))=0.$ 
To complete the proof it remains to apply Proposition \plref{p3.6}.
      \end{proof}
Similarly one proves the following 
   \begin{prop}\label{pr5.38A}
Let $A$, $\ch$ and $R$ be as in Proposition \plref{pr5.32A}. 
Then under the conditions of Proposition \plref{pr5.38} the formal 
deficiency indices of the equation \eqref{5.60} satisfy the inequality  
$\cn_{\pm}(P_+) \le 2n-1.$
    \end{prop}
      \begin{cor}\label{cor5.36} 
Let  $0<c_1\le A(x)\le c_2$\ \ and\ \ $\int_1^{\infty }x^2||R(x)||dx<\infty.$
If 
  \begin{equation}\label{5.62A}
\int_0^\infty (\tr \ch(x))^{1/2}dx=\infty, 
   \end{equation}
then $\cn_{\pm}(P_+) \le 2n-1.$
     \end{cor}
    \begin{proof}
Applying the Cauchy-Bunyakovskii inequality one gets 
      \begin{equation}\label{5.63A}
\int_1^\infty (\tr \ch(x))^{1/2}dx=\int_1^\infty (x^2\tr\ch(x))^{1/2}\cdot 
\frac1x dx\le \int_1^\infty x^2\tr\ch(x)dx.
      \end{equation}
Combining \eqref{5.62A} with \eqref{5.63A} and taking the obvious
inequality $c_1x\cdot I_n\le {\tilde A}(x)\le c_2x\cdot I_n$  into
account one gets 
$\int^{\infty}_a\tr\bigl({\widetilde A}(x)H(x){\widetilde A}(x)\bigr)dx=\infty$

To complete the proof it remains to apply Proposition \plref{pr5.32A}
and note that $\cn_+(P_+)=\cn_-(P_+)$ (see the proof of 
Proposition \plref{pr5.38}).   
    \end{proof}
    \begin{cor}\label{ML-S5.19}
Consider the scalar $(n=1)$ equation \eqref{5.60}. 
Let  $A$ and $R$ satisfy  the conditions \eqref{5.60A}.
Then:
   \begin{thmenum}
\item $\cn_{\pm }(P_+)=1$ if and only if $\int ^{\infty }_0
({\widetilde A}^2(x)+1)\ch(x)dx = \infty$.
\item $\cn_{\pm }(P_+)=2$ if and only if  $\int ^{\infty }_0
({\widetilde A}^2(x)+1)\ch(x)dx < \infty$.
     \end{thmenum}
        \end{cor}
   \begin{proof}
1) By \eqref{3.2a} $\cn_{\pm}(P_+)\ge 1.$ On the other hand by Proposition
\plref{pr5.32A} either $\cn_+(P_+)< 2$  or $\cn_-(P_+) < 2.$ Since 
maximum values of the formal deficiency indices are  attained only 
simultaneosly, one gets  $\cn_{\pm}(P_+)=1.$

2)This assertion is a special case of Proposition \plref{pr5.32A}.
    \end{proof}

  \begin{remark}\label{re5.40} 
Consider the scalar equation \eqref{5.60}. In {\cite{Kre:TFO}} 
(see also {\cite{KacKre:SFS}})
M. Krein stated (without proof) the following result: 

if $n=1, A=1, R$ is semibounded below and 
$\int_1^{\infty }\ch (x)^{1/2}dx=\infty $ then $\cn_{\pm}(P_+)=1.$ 

It follows from Proposition \plref{pr5.33A} that this result 
fails. Moreover, we have explicit counterexamples:
  \begin{equation*}  
R=-k^2<0, \quad \ch(x)=(1+x)^{-2+\eps} \   (0\le \eps <1).
  \end{equation*}
In this case by Proposition \plref{pr5.33A} $\cn_{\pm }(P_+)=2$, but 
$\int_1^{\infty}\sqrt{\ch (x)}dx=\infty .$

Nevertheless Krein's result remains valid for $R\ge 0$ 
(see Remark \plref{r5.04}).
We emphasize however that the statements of Propositions \plref{ML-S5.16}, 
\plref{pr5.32A} and \plref{pr5.33A}  are stronger than the statements we 
obtain by applying Theorem \plref{th5.4} to \eqref{5.35} and \eqref{5.60} 
respectively.
In particular, for $n=1$ these statements are stronger than Krein's result. 
Say, if in Corollary \plref{ML-S5.19} $A=1,\ \ch (x)=(1+x)^{-3}$ then 
$\int_1^{\infty}x^2\ch (x)dx=\infty $ and $\cn_{\pm }(P_+)=1$, but
$\int_1^{\infty }\ch (x)^{1/2}dx<\infty .$
    \end{remark}

\comment{
{123}}

\bibliography{mlabbr,books99,papers99}
\bibliographystyle{lesch}

\end{document}